\documentclass[10pt]{amsart}
\usepackage{color}
\usepackage{bbm}
\usepackage{mathtools}

\allowdisplaybreaks[4]

\definecolor{c20}{rgb}{0.,0.7,0.}
\definecolor{c30}{rgb}{0.,0.,1.}
\definecolor{c40}{rgb}{0.3,0.3,0.9}
\definecolor{c50}{rgb}{1,0,0}
\definecolor{c60}{rgb}{1,0.9,0.1}

\def\kk#1{\textcolor{cyan}{#1}}
\def\kk#1{#1}
\def\qq#1{\textcolor{magenta}{#1}}
\def\qq#1{#1}

\def\ccP#1{\textcolor{c30}{#1}}
\def\ccP#1{#1}

\def\rJ#1{\textcolor{c30}{#1}}
\def\rJ#1{#1}

\def\kls#1{\textcolor{red}{#1}}
\def\ksn#1{\textcolor{red}{#1}}
\def\kls#1{#1}
\def\ksn#1{#1}

\newcommand{\kb}[1]{\boldsymbol{#1}}
\newcommand{\vk}[1]{\kb{#1}}

\newcommand{\abs}[1]{\left\lvert #1 \right\rvert}
\newcommand{\Abs}[1]{ \biggl \lvert #1 \biggr \rvert}
\newcommand{\ABs}[1]{ \biggl \lvert #1 \biggr \rvert}

\newcommand{\E}[1]{\mathbb{E}\left(#1\right)}
\newcommand{\VA}[1]{\mathbb{V}ar\left(#1\right)}

\newcommand{\pk}[1]{\mathbb{P} \left( #1 \right ) }

\newcommand{\R}{\mathbb{R}}

\newcommand{\N}{\mathbb{N}}

\newcommand{\inn}{\in \N}
\newcommand{\ldot}{,\ldots,}

\newcommand{\BQN}{\begin{eqnarray}}
\newcommand{\EQN}{\end{eqnarray}}
\newcommand{\BQNY}{\begin{eqnarray*}}
\newcommand{\EQNY}{\end{eqnarray*}}

\newcommand{\BS}{\begin{sat}}
\newcommand{\ES}{\end{sat}}
\newcommand{\BT}{\begin{theo}}
\newcommand{\ET}{\end{theo}}
\newcommand{\BK}{\begin{korr}}
\newcommand{\EK}{\end{korr}}

\newcommand{\BD}{\begin{de}}
\newcommand{\ED}{\end{de}}

\newcommand{\BIT}{\begin{itemize}}
\newcommand{\EIT}{\end{itemize}}
\newcommand{\BDI}{\begin{description}}
\newcommand{\EDI}{\end{description}}

\newcommand{\BRM}{\begin{remarks}}
\newcommand{\ERM}{\end{remarks}}

\newcommand{\BEL}{\begin{lem}}
\newcommand{\EEL}{\end{lem}}

\newcommand{\BES}{\begin{sat}}
\newcommand{\EES}{\end{sat}}

\newtheorem{theo}{Theorem}[section]
\newtheorem{sat}[theo]{Proposition}
\newtheorem{de}[theo]{Definition}
\newtheorem{lem}[theo]{Lemma}

\newtheorem{example}[theo]{Example}
\newtheorem{korr}[theo]{Corollary}
\newtheorem{remark}[theo]{Remark}
\newtheorem{remarks}[theo]{Remarks}

\newcommand{\COM}[1]{}

\newcommand{\QED}{\hfill $\Box$}

\def\IF{\infty}

\def\LT{\left}
\def\RT{\right}
\def\H{\mathcal{H}}

\def\ooo{(1+o(1))}

\def\H{\mathcal{H}}

\def\vn{\varepsilon}
\def\Var{\text{Var}}

\def\Del{\triangle}

\def\wtb{\widetilde{\vk b}}
\def\olv{\overline{\vk v}}

\def\mcalC{\mathcal C}
\def\hDel{\widehat\Delta}
\def\hN{\widehat N}
\def\hp{\widehat p}
\def\hpi{\widehat \pi}

\begin{document}

\title[\kls{Probability of entering the upper orthant}]
{Probability of entering an orthant by correlated fractional Brownian motion
with drift: exact asymptotics}

\author{Krzysztof D\c{e}bicki}
\address{Krzysztof D\c{e}bicki, Mathematical Institute, University of Wroc\l aw, pl. Grunwaldzki 2/4, 50-384 Wroc\l aw, Poland}
\email{Krzysztof.Debicki@math.uni.wroc.pl}

 \author{Lanpeng Ji}
\address{Lanpeng Ji, School of Mathematics, University of Leeds, Woodhouse Lane, Leeds LS2 9JT, United Kingdom
}
\email{l.ji@leeds.ac.uk {\it (This is the corresponding author.)}}

\author{Svyatoslav Novikov}
\address{Svyatoslav Novikov, Department of Actuarial Science, 
University of Lausanne, UNIL-Dorigny, 1015 Lausanne, Switzerland
}
\email{Svyatoslav.Novikov@unil.ch}

\bigskip

\date{\today}
 \maketitle

{\bf Abstract:}
For  $\{\vk{B}_{H}(t)= (B_{H,1}(t) \ldot B_{H,d}(t))^{\kls{\top}},t\ge0\}$,
where
$\{B_{H,i}(t),t\ge 0\}, 1\le i\le d$ are mutually independent fractional Brownian motions,
we obtain the exact asymptotics of
$$
\mathbb P (\exists t\ge 0: A \vk{B}_{H}(t) - \vk \mu t >\vk \nu u), \ \ \ \ u\to\infty,
$$
where $A$
is a  non-singular $d\times d$ matrix and
$\vk\mu=(\mu_1,\ldots, \mu_d)^{\kls{\top}}\in \R^d$,
$\vk\nu=(\nu_1, \ldots, \nu_d)^{\kls{\top}} \in \R^d$ are such that  there exists some $1\le i\le d$ such that
$\mu_i>0, \nu_i>0.$

\medskip


\medskip


\section{Introduction}

Consider a vector-valued Gaussian  process
$\{{\vk{X}}(t),t\ge 0\}$, where ${\vk{X}}(t)=A \vk{B}_{H}(t)$ with
$A\in \R^{d\times d}$ a non-singular matrix and  $\{\vk{B}_{H}(t)= (B_{H,1}(t) \ldot B_{H,d}(t))^{\kls{\top}},t\ge0\}$
\ccP{with} $\{B_{H,i}(t),t\ge 0\}, 1\le i\le d$ \qq{($d\inn$)} \ccP{being} mutually independent fractional Brownian motions (fBms),
i.e., centered Gaussian processes with stationary increments, continuous sample paths and variance functions
$\Var(B_{H,i}(t))=t^{2H}$, $H\ccP{\in(0,1)}$.

We focus on the exact asymptotic behavior of the
probability
that a drifted correlated fractional Brownian motion $\vk{X}$ enters orthant
$\mathcal{O}_u=\{(x_1 \ldot x_d)^{\kls{\top}}: x_i> \nu_i u, i=1 \ldot d\}$ over an infinite-time horizon, i.e.,
\begin{eqnarray}
P(u):=\pk{ \exists t\ge 0:  \vk{X}(t) - \vk \mu t \ \qq{\in \mathcal{O}_u}}=\pk{\exists_{t \ge 0}
 \forall_{i=1 \ldot d} X_{i}(t) -\mu_i t> \nu_iu },
 \label{aim}
\end{eqnarray}
as $u\to\infty$
for  $\vk\mu=(\mu_1,\ldots, \mu_d)^{\kls{\top}}\in \R^d$ and $\vk\nu=(\nu_1, \ldots, \nu_d) ^{\kls{\top}}\in \R^d$.

We are interested in the case that the above probability is a rare event,
that is, $P(u)\to 0$ as $u\to\IF$, for which we shall assume \qq{that} there exists some $1\le i\le d$ such that
$
 \mu_i>0, \nu_i>0.
$

The probability $P(u)$ defined in \eqref{aim}
is of interest both for theory-oriented studies and for applied-mathematics problems. 
One of important motivations to analyze \eqref{aim} stems from {\it ruin theory}, where
$P(u)$ describes simultaneous ruin probability in infinite-time horizon of $d$ dependent business lines
whose risk processes $R_i(t),t\ge 0$ are modeled by
\[
R_i(t)=\nu_iu  +\mu_i t -X_{i}(t),
\]
where $\nu_iu $ represents the initial capital, $\mu_i$ is the \qq{net profit rate} and
$X_{i}(t) $ is the \qq{net loss up to time $t$}; we refer to \cite{Mic99}
for the formal justification of the use of fractional Brownian motion to model the risk process.

In the 1-dimensional case, $d=1$, the exact asymptotics for $P(u)$
was derived in the seminal paper by H\"{u}sler \& Piterbarg \cite{HP99}; see also \cite{Deb99, HuP04, Die05}
for extensions to other classes of \ccP{stochastic}
processes with stationary increments.

In the multidimensional case,  the exact asymptotics of $P(u)$ as $u\to\infty$ is known only for the special Brownian model, i.e., when $H=1/2$; see \cite{DHJRBM}.  The strategy of the proof there, although in its roots based on {\it the double sum} technique developed in the 1-dimensional setting for extremes of Gaussian processes (see, e.g. \cite{PickandsA, PicandsB,Pit96}),
needed new ideas that in several key steps of the argumentation significantly differ from methods used in the 1-dimensional case.
In particular, one of the difficulties is the lack of Slepian-type inequalities that could be applied in the multidimensional setting, which was overcome  in \cite{DHJRBM} by the heavy use of the independence of increments property of Brownian motion.

In this contribution, we aim to complement the findings of \cite{DHJRBM} by tackling the fBm problem \eqref{aim}
for $H\neq 1/2$.  Interestingly, in contrary to the Brownian case, the full analysis of all the cases
needs to consider two separate scenarios
described by the local behaviour of the function
\[
\vk{\zeta}(t)=\ccP{(\zeta_1(t), \ldots, \zeta_d(t))^{\kls{\top}}} := (\vk{\nu}+\vk{\mu}t)/t^H, t>0
\]
in the neighbourhood of the unique point $t_0$ that minimizes function
\begin{eqnarray}
g(t):=\frac{1}{t^{2H}} \inf_{\vk{v} \ge \vk\nu+\vk{\mu} t}  \vk{v}^\top (AA^\top)^{-1}  \vk{v} \label{g1}
\end{eqnarray}
over $t\ge0$,  where the point \qq{$t_0u$} has a natural interpretation as the {\it most probable time}
for the process $\{\vk{X}(t) - \vk \mu t,t\ge0\}$ to enter orthant $\mathcal{O}_u$;
see also Section \ref{s.prem}. Let $I\subset\{1,...,d\}$ be the set of coordinates that contribute to the
asymptotics of $P(u)$, as $u\to\infty$; see Section \ref{s.prem} for the details \qq{about} how to specify $I$.\\
In the first scenario, when $H<{1}/{2}$, or $H>1/2$ and $\ccP{\zeta}^{'}_i(t_0)={0}$ for all $i\in I$,
the local steepness of \qq{the correlation function} of $\ccP{X}_i$ is higher than
the local steepness of  $\ccP{\zeta}_i(t)$ in the neighbourhood of the point
$t_0$. This case can be solved by an adaptation of the technique developed in
\cite{DHW} for extremes of centered vector-valued Gaussian processes over a finite time horizon.\\
The complementary case, $H>1/2$ but  $\ccP{\zeta}^{'}_i(t_0)\neq{0}$ for some $i\in I$, is different from the previous one
since for those coordinates $i\in I$ for which  $\ccP{\zeta}^{'}_i(t_0)\neq{0}$
the local steepness of \qq{the correlation function of} $\ccP{X}_i$ is lower than
the local steepness of $\ccP{\zeta}_i(t)$  in the neighbourhood of the point
$t_0$.
This scenario needs a novel approach which leads also to a different asymptotics of $P(u)$ as $u\to\infty$.

\qq{The results derived in this contribution go in line with recent findings on the tail asymptotics of extremes of
vector-valued Gaussian processes, where most of the available literature deals with centered marginal processes
or over a compact parameter space \cite{DHW, BDK23,Iev24, DIK23,ChX23}.
}
\medskip

{\bf Notation.} We shall use some standard notation which is common when dealing with (column) vectors.
All the operations on vectors are meant componentwise, for instance, for any given
 $ \vk{x} = (x_1,\ldots,x_d)^{\kls{\top}}\in \R ^d$ and $\vk{y} = (y_1,\ldots,y_d)^{\kls{\top}} \in \R ^d $, we write $ \vk{x} > \vk{y} $ if and only if
  $ x_i > y_i $  for all $ 1 \leq i \leq d $, {write $1/\vk{x}=(1/x_1,\cdots,1/x_d)^{\kls{\top}}$ if $x_i\neq 0, 1 \leq i \leq d$}, and write $\vk{x}\vk{y}=(x_1y_1 \ldot x_dy_d)^{\kls{\top}}$
  and $a\vk x =(a x_1,\ldots, ax_d)^{\kls{\top}},$ $a\in \R.$ 
Further, we set  $ \vk{0}  := (0,\ldots,0)^{\kls{\top}} $ and $ \vk{1} : = (1,\ldots,1)^{\kls{\top}}$ whose dimension will be clear from the context. Moreover, denote $\abs{\vk x}$ as the $L_1$-norm of  $ \vk{x} = (x_1,\ldots,x_d)^{\kls{\top}}\in \R ^d$.

If $I\subset \{1,\ldots,d\}$, then for a vector $\vk a\in\R^d$ we denote by $\vk{a}_I=(a_i, i\in I)$ a sub-block vector of $\vk a$. Similarly, if further $J \subset \{1\ldot d\}$, for  a matrix $M=(m_{ij})_{i,j\in \{1,\ldots,d\}}\in \R^{d\times d}$ we denote by $M_{IJ}=(m_{ij})_{i\in I, j\in J}$ the sub-block matrix of $M$
determined by $I$ and $J$. Further, write $M_{II}^{-1}=(M_{II})^{-1}$ for the inverse matrix of $M_{II}$ whenever it exists.
Denote by $\abs{I}$ the number of elements in the index set $I$ and by $\abs{M}$ the determinant of a square matrix $M$.

For two positive functions $f, h$ and some $u_0\in[-\IF, \IF]$, write  $ f(u) = h(u)(1+o(1))$ or $h(u)\sim  f(u)$ if $ \lim_{u \to u_0} f(u) /h(u)  = 1 $,  write $ f(u) = o(h(u)) $ if $ \lim_{u \to u_0} {f(u)}/{h(u)} = 0$. 
\kls{Moreover, for an event $E$ we denote $\mathbb{I}_{(E)}$ as the indicator function of $E$.}


\medskip

{\bf Organization of the paper.}
Some preliminary results related with properties and the role of function $g$ defined in \eqref{g1} are presented in Section 2. The main result of this contribution, which is Theorem \ref{Thm:main}, is given in
Section 3, followed by an illustrative example.
In Section 4, we give the proof of the main result.
All \ccP{other} technical proofs are relegated to \ccP{Appendix}.

\section{ Preliminary results}\label{s.prem}

%


\COM{ 
We shall refer to $\{\vk{X}(t), t\ge0\}$
as
a centered $n$-dimensional {\it vector-valued} Gaussian process,
where $ \vk{X}(t)=(X_1(t) \ldot X_n(t))$ with $X_i$'s being \ccP{independent} centered Gaussian processes with a.s. continuous sample paths.
Since $n$ hereafter is always fixed we shall occasionally omit
"$n$-dimensional", mentioning simply that $\vk{X}$ is a centered vector-valued Gaussian \rrJ{process}.
Define next
\begin{eqnarray*}
\sigma_{\vk{X}}^2(\cdot)   =   (\sigma_{X_1}^2(\cdot),\ldots,\sigma_{X_n}^2(\cdot)), \quad  R_{\vk{X}}(\cdot,\cdot)  =  (R_{X_1}(\cdot,\cdot),\cdots,R_{X_n}(\cdot,\cdot)),
\end{eqnarray*}
with $\sigma_{X_i}^2(t)  =  \mathrm{Var}(X_i(t)) $ and $R_{X_i}(s,t)  =  \mathrm{Cov}(X_i(s),X_i(t))$.

Let in the following
$\{B_{i,\kappa}(t),t\in \R\}, 1\le i\le n$  be $n$ mutually independent standard fractional Brownian motions (fBm's) defined on $\R$ with \ccP{common} Hurst
index $\kappa/2\in(0,1]$, and set $\vk{B}_\kappa(t)= ( B_{1, \kappa}(t) \ldot B_{n,\kappa}(t)).$

A key step in the investigation of the tail asymptotics of supremum of Gaussian processes is the derivation of the tail asymptotic behaviour of the supremum taken over
"short intervals". For the stationary case this is achieved by the so-called Pickands lemma.
The non-stationary case is covered by the so-called Piterbarg lemma
(see \cite{DHL14Ann, DHJParisian, HJ14d} for similar terminology and related results).
Before deriving an extension of these classical results for the vector-valued \ccP{Gaussian} processes, we need to introduce some
further notation.

} 

It is known that the approximation of the probability $P(u)$
depends on the solution of a related quadratic optimization problem.
In particular,
in the light of \cite[Theorem 1]{RolskiSPA}, the logarithmic
asymptotics   can be derived and takes the following form 
\def\gt{\widehat g}
\BQN \label{eq1D}
-\ln P(u) \sim  \frac{\gt}{2} u^{2(1-H)}, \quad \text{with}\ \  \gt=\inf_{t\ge 0} g(t),
 \EQN
 where
\BQN \label{fgt}
g(t)=\frac{1}{t^{2H}} \inf_{\vk{v} \ge \vk\nu+\vk{\mu} t}  \vk{v}^\top \Sigma^{-1}  \vk{v}, \quad
\Sigma=AA^\top.
\EQN

The properties of function $g$, in particular, the existence of its \kk{minimizer} and expansions \kk{in the neighbourhood of} this point (when exists) are crucial to the exact asymptotic analysis.
 In order to introduce some further notation and for further reference, we present a lemma on the quadratic optimization problem stated in  \cite{HA2005} (see also  \cite{ENJH02}). 

\BEL \label{AL}
Let
$\Sigma \in \R^{d \times d},d\ge 2$ be a positive definite  matrix. 
If  $\vk{b}\in  \R ^d \setminus (-\infty, 0]^d $, then the quadratic programming problem
$$ P_\Sigma(\vk{b}): \text{minimize $ \vk{x}^\top \Sigma^{-1} \vk{x} $ under the linear constraint } \vk{x} \ge \vk{b} $$
has a unique solution $\widetilde{\vk{b}}$ and there exists a unique non-empty
index set $I\subseteq \{1, \ldots, d\}$ such that
\BQN  \label{eq:IJi}
&&\widetilde{\vk{b}}_{I}=
\vk{b}_{I}\not= \vk 0_I, \quad \widetilde {\vk{b}}_{I^c} = \Sigma_{I^cI}\Sigma_{II}^{-1} \vk{b}_{I}\ge \vk{b}_{I^c}, \quad \Sigma_{II}^{-1} \vk{b}_{I}>\vk{0}_I,\\
&&\min_{\vk{x} \ge  \vk{b}}
\vk{x}^\top \Sigma^{-1}\vk{x}= \widetilde{\vk{b}}^\top \Sigma^{-1} \widetilde {\vk{b}}   =
\vk{b}_{I}^\top \Sigma_{II}^{-1}\vk{b}_{I}>0,  
 \label{eq:new}
\EQN
where $I^c= \{1, \ldots, d\} \setminus I.$
Moreover, denoting $\vk w=\Sigma^{-1} \widetilde{\vk{b}}$ we have
\BQN\label{eq:ww}
\vk w_I=\Sigma_{II}^{-1} \vk{b}_{I}>\vk{0}_I, \quad \vk w_{I^c} =\vk 0_{I^c}.
\EQN
\EEL

The next lemma includes some properties of the function $g$ and its relative $g_I(t)=\frac{1}{t^{2H}}  (\vk{\nu}+\vk \mu t)^\top_{I} \Sigma_{II}^{-1}   (\vk{\nu}+\vk \mu t)_I$ with $I$
\kk{the index set as in Lemma \ref{AL}.}
\kk{We defer its proof to Appendix.}

\BEL \label{lem:gc1}
\kk{Function  $g\in C^1(0,\IF)$ and}  achieves its unique minimum at
\BQN\label{eq:t0}
t_{0}=\frac{\sqrt{4 (\vk\nu_{I}^\top \Sigma^{-1}_{II}  \vk{\mu}_{I} )^2(1-2H)^2 +16 H(1-H) \vk\nu_{I}^\top \Sigma^{-1}_{II}  \vk\nu_{I} \vk{\mu}_{I}^\top \Sigma^{-1}_{II}  \vk{\mu}_{I} }- 2 (1-2H)  \vk\nu_{I}^\top \Sigma^{-1}_{II}  \vk{\mu}_{I}}{4 \vk{\mu}_{I}^\top \Sigma^{-1}_{II}  \vk{\mu}_{I} (1-H)} >0  
\EQN
\kk{with} 
\BQN \label{eq:intr1}
g(t_0)=\inf_{t>0}\frac{1}{t^{2H}} \inf_{\vk{v} \ge \vk\nu+\vk{\mu} t}  \vk{v}^\top \Sigma^{-1}  \vk{v}=\frac{1}{t_0^{2H}}  \vk{b}^\top_{I} \Sigma_{II}^{-1}  \vk{b}_{I}=g_I(t_0),
\EQN
where
$$\vk{b}=\vk b(t_0), \quad \text{with \ } \vk b(t)= \vk\nu + \vk \mu t,$$
 and $I=I(t_0)$ being the index set corresponding to the solution of $P_\Sigma (\vk b)$.  Moreover,
   \BQN
g_I(t_0+ t) = g_I(t_0)+ \frac{g''_I(t_0)}{2} t^2\ooo, \ \ \ \ t\to 0, \label{eq:gt0pm}
\EQN
where
$$g''_I(t_0)=\frac{1}{t_0^{2H+1}}\LT( 4 \vk\mu_{I}^\top \Sigma^{-1}_{II}  \vk\mu_{I} (1-H) t_0 +2(1-2H)\vk \nu_{I}^\top \Sigma^{-1}_{II}  \vk\mu_{I}\RT) >0.$$
\EEL

\begin{remark}
 Note that $t_0$ given in \eqref{eq:t0} is actually an equation of $t_0$ because $I=I(t_0)$ is a set function of $t_0$. Here, for any fixed $t>0,$  $I(t)\subseteq \{1,2,\ldots, d\}$ is the index set of the solution to the quadratic programming problem $P_{\Sigma}(\vk b(t))$, see Lemma \ref{AL}. We remark that, for specific problems,
both the index set $I$ and $t_0$ can be identified explicitly; see \kk{Example \ref{Ex:1} in Section \ref{s.main}}
or the examples presented in \cite{DHJRBM}.
%
\end{remark}

Hereafter we shall use the notation $\vk{b}=\vk{b}(t_0)= \vk\nu + \vk \mu t_0$, and $I=I(t_0)$ for the {\it essential} index set of the quadratic programming problem $P_{\Sigma}(\vk b)$.  Furthermore,
let  $\widetilde{\vk{b}}$ be the unique solution
of $P_{\Sigma}(\vk b)$. 
If $I^c=\{1,\ldots,d\}\setminus I\neq\emptyset $, we define
the {\it weakly essential index} and the {\it unessential index} sets by
\BQN\label{K.def}
K=\{j\in I^c: \widetilde{b}_j=\Sigma_{jI}\Sigma_{II}^{-1}\vk{b}_I=b_j\},
 \quad
\text{  and } J=\{j\in I^c: \widetilde{b}_j=\Sigma_{jI}\Sigma_{II}^{-1}\vk{b}_I>b_j\},
\EQN
respectively.


\section{Main Result }\label{s.main}

\kk{In this section we present the main result of this contribution.}
\qq{Recall that through the whole paper we assume \qq{that} there exists some $1\le i\le d$ such that
$ \mu_i>0, \nu_i>0.$}
Denote $\vk W_{I}(t) =D \vk B_{H,I}(t)$, with $D$ a  matrix such that $D D^\top = \Sigma_{II}$.
We define a {\it generalized Pickands constant} as
\BQNY
\mathcal H_{I}= \lim_{T\to \IF} \frac{1}{T} \mathcal H_{I}(T) \in(0,\IF),
\EQNY where
\BQNY
\mathcal H_{I}(T)= \int_{\R^{\abs{I}}} e^{\frac{1}{t_0^{2H}} \vk w_I^\top \vk x_I}   \pk{\exists_{t \in[0, T]}
\vk W_{I}(t)  -  \frac{1}{2t_0^{2H}}\vk b_I t^{2H}>\vk x_I } d\vk x_I,  
\EQNY
\rJ{with $\vk b = \vk\nu + \vk \mu t_0$  and $\vk w_I =\Sigma_{II}^{-1} \vk{b}_I.$
We remark that $\mathcal{H}_I$ is well-defined, finite and positive, since it is a multiple of the multidimensional Pickands constant $\mathcal{H}_{2H,V}$ defined in (2.5) of \cite{DHW} with $V=\LT(2 t_0^{4H}\RT)^{-1} \mathrm{diag}(\vk{w}_I) \Sigma_{II} \mathrm{diag}(\vk{w}_I)$. }
%
%

For $K$ defined in \eqref{K.def}, let 
\BQN\label{eq:mcalC}
&& \kls{\mcalC_{K}}:= 
\left\{\begin{array}{lll}
\kls{\frac{1}{\sqrt{ (2\pi  t_0^{2H} )^{\abs{I}}  \abs{\Sigma_{II}}}} }
\int_{\R} e^{-\frac{1}{4} g_I''(t_0) y^2} \ \cdot \pk{\vk Y_K< t_0^{-H} \LT(\vk \mu_K - \Sigma_{KI}\Sigma_{II}^{-1} \vk \mu_I\RT) y} dy,
& K\neq \emptyset,\\
\kls{\sqrt{ \frac{4\pi}{g''_I(t_0) (2\pi  t_0^{2H} )^{\abs{I}}  \abs{\Sigma_{II}}} } },
&  K=\emptyset.
\end{array}\right.
\EQN
\kls{where  $\vk Y_K \overset{d}\sim \mathcal N\LT(\vk 0_K, \Sigma_{KK} -  \Sigma_{KI}\Sigma_{II}^{-1} \Sigma_{IK} \RT)$}.

\qq{
Following Section \ref{s.prem}, the logarithmic asymptotics of $P(u)$
as $u\to\infty$ depends on $t_0$, the minimizer of function $g$
defined in \eqref{fgt}.
As stated in the following theorem, 
the exact asymptotics of $P(u)$
splits on two scenarios.
In order to catch an intuition behind this division, for a while let us consider
the 1-dimensional 
problem
	$
	\pk{\exists_{t \ge 0} B_{H,i}(t) -\mu_i t>\nu_i u},
	$
\ccP{assuming $\mu_i,\nu_i>0.$}
By self-similarity of $B_{H,i}$ we get
\begin{eqnarray*}
\pk{\exists_{t \ge 0} B_{H,i}(t) -\mu_i t>\nu_i u}
=
\pk{\exists_{t \ge 0} \frac{B_{H,i}(t)}{\nu_i u +\mu_i t}>1}
=
\pk{\exists_{t \ge 0} \frac{B_{H,i}(t)}{\nu_i +\mu_i t}>u^{1-H}}
\end{eqnarray*}
and thus, following the same lines of reasoning as in Section \ref{s.prem} but for 1-dimensional setting,
the logarithmic asymptotics of the above is determined by $t_{0,i}$, the unique minimizer of
$\zeta_i(t)=\frac{\nu_i +\mu_i t}{t^H}$, $t>0$, that is the point that satisfies
${\zeta_i}'({t_{0,i}}) = {0}$ or equivalently
\BQN\label{eq:Hmui}
H {\nu_i} = (1-H) {t_{0,i}} {\mu_i}.\EQN
This leads to two cases, where the play between the value of the optimizing point $t_0$
and the optimizers $t_{0,i}$ of
$\zeta_i(t)$ 
for $i\in I$
is crucial: \\
$\diamond$ $H<{1}/{2}$, or $H>1/2$ and $\ccP{\zeta}^{'}_i(t_0)={0}$ for all $i\in I$.
Then the local steepness of \qq{the correlation function} of $\ccP{X}_i$ is higher than
the local steepness of  $\ccP{\zeta}_i(t)$ in the neighbourhood of the point
$t_0$.\\
$\diamond$ $H>1/2$ but  $\ccP{\zeta}^{'}_i(t_0)\neq{0}$ for some $i\in I$.
For the coordinates $i\in I$ for which  $\ccP{\zeta}^{'}_i(t_0)\neq{0}$
the local steepness of \qq{the correlation function of} $\ccP{X}_i$ is lower than
the local steepness of $\ccP{\zeta}_i(t)$  in the neighbourhood of the point
$t_0$.
}

\kk{The following theorem constitutes the main finding of this contribution.}
\ccP{Let us recall notation $x_{-}=\max(0,-x),\,x \in \R$. }

\BT \label{Thm:main}
\begin{itemize}
\item[(i).] If $H<{1}/{2}$, or $H>1/2$ and $H \vk{\nu}_{I} = (1-H) t_0 \vk{\mu}_{I}$, then
\BQNY
P(u) \sim    \mathcal H_{I}\ksn{\mcalC_{K}} u^{-\abs{I}(1-H)+1/H+H-2} e^{-\frac{g_I(t_0)}{2}u^{2(1-H)} },\ \ u\to\IF.
\EQNY
\item[(ii).] If $H>1/2$ but $H\vk{\nu}_{I} \neq (1-H) t_0\vk{\mu}_{I}$, then
\BQNY
P(u) \sim
\ksn{\frac{t_0^{2H(\abs{I}-1)}}
{\prod_{i\in I}w_i}
  \mcalC_{K}} \LT(\sum_{i \in I} \left(w_i \mu_i - \frac{H}{t_0} w_i b_i\right)_{-} \RT)  u^{-\abs{I}(1-H)+1-H}
 e^{-\frac{g_I(t_0)}{2}u^{2(1-H)} }, \ \ u\to\IF,
\EQNY
where
$$
\sum_{i \in I} \left(w_i \mu_i - \frac{H}{t_0} w_i b_i\right)_{-} >0.
$$
\end{itemize}

\ET
\begin{remark}
 We remark on the role of the index set $J$ played in the asymptotics when it is non-empty. It is concluded from \cite{DHJRBM} that  the index set $J$ and the corresponding components do not play any role in the exact asymptotics of $P(u)$ for $H=1/2$. 
 \kls{It follows from Theorem \ref{Thm:main} that the same observation applies also for $H \neq 1/2$.}
\COM{
  Indeed, this can be seen by inserting the second scenario value of \eqref{eq:mcalC} and noting that
$
\abs{\Sigma_{II}} = \abs{\Sigma} \abs{(\Sigma^{-1})_{JJ}} .
$ However, if $K\neq \emptyset$, then the inner integration in the first scenario of \eqref{eq:mcalC} seems difficult to simplify (unless $(\Sigma^{-1})_{I^cI^c}$ is of a simple form, e.g., a diagonal block matrix) and thus it is hard to conclude in general  whether $J$ is not playing any role. Note in passing that if $K\neq \emptyset$ and $J= \emptyset$,  then the inner integration in the first scenario of \eqref{eq:mcalC} becomes
$$
\sqrt{ (2\pi  t_0^{2H} )^{\abs{K}} }/ \sqrt{ \abs{(\Sigma^{-1})_{KK}}} \ \cdot \pk{\vk Y_K< t_0^{-H} \LT(\vk \mu_K - \Sigma_{KI}\Sigma_{II}^{-1} \vk \mu_I\RT) y},
$$
with $\vk Y_K \overset{d}\sim \mathcal N\LT(\vk 0_K, \Sigma_{KK} -  \Sigma_{KI}\Sigma_{II}^{-1} \Sigma_{IK} \RT)$. This  leads to a formulation of $\mcalC_{K,\emptyset}$ that is consistent with the constant involving $K$ in \cite{DHJRBM}.
}
\end{remark}

We conclude this section with an illustrative example, where we will see how the index sets $I, K, J$ and the optimal point $t_0$ are derived and how the different cases may appear. Our purpose of this example is not to be as general as possible, but to be restrictive so that it includes an interesting scenario.

\begin{example}\label{Ex:1}

We consider a $4$-dimensional Gaussian process with independent (positive or negative) drifted fBm components. Precisely, let
$$
d=4,\ \ \Sigma=\mathrm{Id},\ \ \nu_i>0, i=1,2,3,4, \ \ \mu_1, \mu_2>0,\ \ \mu_3, \mu_4<0.
$$
We also assume that
\BQN
\IF=:t'_3>t'_2:=\frac{\nu_3}{\abs{\mu_3}} > \frac{\nu_4}{\abs{\mu_4}}=:t'_1 >t'_0:=0.
\EQN
Denote
$$
I_1=\{1,2,3,4\},  \ \ I_2=\{1,2,3\}, \ \ I_3=\{1,2\}, \ \ I^c_j=\{1,2,3,4\}\setminus I_j,\  j=1,2,3.
$$
It can be seen that
$$
\vk \nu_{I_j} +t \vk \mu_{I_j} >\vk 0_{I_j},\  \  \vk \nu_{I_j^c} +t \vk \mu_{I_j^c} \le \vk 0_{I_j^c}, \ \ \ \  t\in [t'_{j-1}, t'_j)
$$
and thus, by Lemma \ref{AL},
$$
I(t) = I_j, \  t\in [t'_{j-1},t'_j), \ \ \ \ j=1,2,3.
$$
Further, it follows from Lemma \ref{lem:gc1} that the optimal point $t_0$ is equal to the
$t_0^{(k)}$ defined as \eqref{eq:t0} with $I_k$, such that (see \eqref{eq:intr1})
$$
g_{I_k}(t_0^{(k)})= \min_{j=1,2,3} g_{I_j}(t_0^{(j)}).
$$
For illustration purpose, we shall assume that we have chosen the model parameters such that $k=3$. Thus, the essential index set is given by $I=I_3=\{1,2\}$ and
$$
t_0=t_0^{(3)} =\frac{\sqrt{4(\sum_{i\in I} \nu_i\mu_i)^2(1-2H)^2+16 H(1-H) \sum_{i\in I} \nu_i^2 \sum_{i\in I} \mu_i^2 }-2(1-2H)\sum_{i\in I} \nu_i\mu_i}{4\sum_{i\in I} \mu_i^2(1-H)}\in[t'_2, \IF).
$$
We further assume that the model parameters were chosen such that $t_0=t'_2$. In such a case, we have
$$
I=\{1,2\}, \ \ K=\{3\}, \ \ J=\{4\}.
$$
Next, let us discuss  different cases distinguished according to $H \vk{\nu}_{I} = (1-H) t_0 \vk{\mu}_{I}$ being valid or not.
\qq{Following the notation introduced at the beginning of this section,} $H \vk{\nu}_{I} = (1-H) t_0 \vk{\mu}_{I}$ means that $t_0=t_{0,1}=t_{0,2}$ (see \eqref{eq:Hmui}) and $\zeta_1'(t_0)=\zeta_2'(t_0)=0.$ In contrast, $H \vk{\nu}_{I} \neq (1-H) t_0 \vk{\mu}_{I}$ means that $t_0$ falls between $ t_{0,1}$ and $t_{0,2}$, and one of $\zeta'_1(t_0)$ and $\zeta'_2(t_0)$ is positive while the other is negative.
Consequently, we obtain the exact asymptotics for
$$\pk{\exists_{t \ge 0}
 \forall_{i=1 \ldot 4} B_{H,i}(t) -\mu_i t> \nu_iu },\ \ u\to\infty,
$$
by applying Theorem \ref{Thm:main}, \ccP{where} 
$$
\ksn{\mathcal C_K} =\ksn{\frac{1}{2\pi t_0^{2H}}} \int_{\R} e^{-\frac{1}{4} g''_I(t_0) y^2} \Phi\LT(\frac{\mu_3}{t_0^H} y\RT) dy,
$$
with  $\Phi(\cdot)$ denoting the \qq{standard normal distribution function}.
\QED

\end{example}

\section{Proof of Theorem \ref{Thm:main}} \label{s.proof}

First note that by self-similarity of the fBms, 
\BQNY
P(u)&=&\pk{\exists_{t \ge 0}
 \vk X(t) -\vk \mu t> \vk \nu u }\\
&=&\pk{\exists_{t \ge 0}
 \vk X(t) > (\vk \nu +\vk \mu t)u^{1-H} }. 
\EQNY
Hereafter, for simplicity we denote $v=u^{1-H}.$  Furthermore, denote
$$
\Delta_v=[t_0- \ln(v)/v, t_0+ \ln(v)/v],\quad \widetilde \Delta_v=[0,\IF)\setminus [t_0- \ln(v)/v, t_0+ \ln(v)/v].
$$
It follows that
\BQN\label{eq:pPpv}
p(v)\le P(u)\le \Pi(v) +p(v),
\EQN
where
$$
p(v)=\pk{\exists_{t \in \Delta_v} \vk X(t) > (\vk \nu +\vk \mu t) v  }, \quad \Pi(v)=\pk{\exists_{t \in\widetilde \Delta_v}  \vk X(t) > (\vk \nu +\vk \mu t) v  }.
$$

The proof consists of two steps. In Step 1, we obtain the asymptotics of $p(v), v\to\IF$. In Step 2, we derive a suitable upper bound for $\Pi(v)$  for all large enough $v$, \kk{which confirms asymptotic negligibility of
$\Pi(v)$ with respect to $p(v)$ as $v\to\IF$}.
The proof is then completed by combining these results. 
Without loss of generality, we shall only consider the most involved case where \kls{$K\neq \emptyset$.} 
Before delving into all the details, \qq{for any $M\in(0, \IF]$ we introduce}
\BQN\label{eq:CKJM}
\ksn{\mcalC_{K,M}}= \kls{ \frac{1}{\sqrt{ (2\pi  t_0^{2H} )^{\abs{I}} \abs{\Sigma_{II}}}  } \int_{-M}^M } e^{-\frac{1}{4} g_I''(t_0) y^2} \ \cdot \pk{\vk Y_K< t_0^{-H} \LT(\vk \mu_K - \Sigma_{KI}\Sigma_{II}^{-1} \vk \mu_I\RT) y} dy.
\EQN
We note that $\ksn{\mcalC_{K}}$, given in \eqref{eq:mcalC} is actually equal to $\ksn{\mcalC_{K,\IF}}$.

\subsection*{\underline{Step 1: Analysis of $p(v)$}} The idea is to split the interval $\Delta_v$ into smaller intervals.
It turns out that we need to distinguish two different cases (i) and (ii) as stated in Theorem \ref{Thm:main}, for case (i) we shall use intervals of the classical Pickands length, but for case (ii) we need to use intervals of a length that is shorter than the Pickands length. These two cases will be discussed seperately below.

\underline{Case (i): \kk{$H<{1}/{2}$, or $H>1/2$ and $H \vk{\nu}_{I} = (1-H) t_0 \vk{\mu}_{I}$.}} Denote, for any fixed  $T>0$ and $v>0$
\BQNY 
\Del_{j;v}=\Del_{j;v}(T)= [t_0+j Tv^{-1/H}, t_0+(j+1) Tv^{-1/H}],\ \ \ -N_v-1\le j\le N_v,
\EQNY
where $N_v=\lfloor T^{-1} \ln(v) v^{1/H-1}\rfloor$
(we denote by $\lfloor\cdot\rfloor$ the floor function). Also denote $N_{v,M} = \lfloor T^{-1} M v^{1/H-1} \rfloor$ for any $M>0$.
By applying the Bonferroni's inequality, we have 
\BQN\label{eq:thetaT}
p_1(v)\ge p(v) \ge p_{2,M}(v)-\pi_M(v),
\EQN
where
\begin{equation*} 
p_1(v)= \sum_{j=-N_v-1}^{N_v}p_{j;v},\ \
p_{2,M}(v) = \sum_{j=-N_{v,M}}^{N_{v,M}}p_{j;v},\ \
\pi_M(v) = \sum_{ -N_{v,M}\le j< l\le  N_{v,M}}p_{j,l;v},
\end{equation*}
with
\begin{equation*} 
p_{j;v}=\pk{\exists_{t\in\Delta_{j;v}} \vk{X}(t)>
(\vk{\alpha}+\vk{\mu}t)v}
\end{equation*}
and
\BQN\label{eq:pjlv}
p_{j,l;v}=\pk{\exists_{t\in \Del_{j;v}}  \vk{X}(t)> (\vk{\alpha}+\vk{\mu} t)v,
 \ \exists_{t\in \Del_{l;v}} \vk{X}(t)> (\vk{\alpha}+\vk{\mu} t)  v  }.
\EQN

Next, we shall deal with the single-sum $p_1(v)$, $p_{2,M}(v)$ and the double-sum $\pi_M(v)$, respectively. For the asymptotics of the single-sum terms, we shall use the following uniform version of a generalized Pickands lemma. The proof of Lemma \ref{lem:Pick1} is displayed in Appendix.

\BEL \label{lem:Pick1} Fix $T>0$. We have, as $v\to\IF,$
\BQNY
&&\pk{\exists_{t \in[t_0+\tau v^{-1/H}, t_0+(\tau+T) v^{-1/H}]}
 \vk X(t) > (\vk \nu +\vk \mu t) v }\\
&&
\sim  v^{-\abs{I}} \mathcal H_{I}(T)   \kls{ \frac{1}{\sqrt{ (2\pi  t_0^{2H} )^{\abs{I}} \abs{\Sigma_{II}}}  }} e^{-\frac{v^2}{2} g_I(t_0+\tau v^{-1/H})}\\
&&
\quad
\ksn{
\times \pk{\vk Y_K< -t_0^{-H} \LT(\vk \mu_K - \Sigma_{KI}\Sigma_{II}^{-1} \vk \mu_I\RT) (\tau v^{\kls{1}-1/H})} 
},
\EQNY
holds uniformly in $\tau$ such that $\abs{\tau}\leq T(N_v+1)$,
where $K$ is the weakly essential index set defined in \eqref{K.def} and \kls{$\vk Y_K$ is given in \eqref{eq:mcalC}}.
\EEL

With  \kk{Lemma} \ref{lem:Pick1}, it is straightforward to check that, as $v\to\IF,$
\BQN\label{eq:p12v}
&&
p_{2,M}(v)\sim
\ksn{
v^{-\abs{I}+1/H-1} \frac{\mathcal H_{I}(T)}{T}\mathcal{C}_{K,M} e^{-\frac{v^2}{2} g_I(t_0)}
},
\EQN
where $\ksn{\mcalC_{K,M}}$ is given in \eqref{eq:CKJM}. Indeed, by  Lemma  \ref{lem:Pick1} and \eqref{eq:gt0pm}, we derive that, as $v\to\IF,$
\BQNY
p_{2,M}(v)
&\sim&
v^{-\abs{I}+1/H-1}
\frac{\mathcal H_{I}(T)}{T}  \kls{ \frac{1}{\sqrt{ (2\pi  t_0^{2H} )^{\abs{I}} \abs{\Sigma_{II}}}  }}  e^{-\frac{v^2}{2} g_I(t_0)}\\
&&\times
\sum_{j=-N_{v,M}}^{N_{v,M}}
(T v^{1-1/H})
e^{-\frac{g_I''(t_0)(jT v^{1-1/H})^2}{4}}\\
&&\qquad \qquad\ksn{
	\times \pk{\vk Y_K< - t_0^{-H} \LT(\vk \mu_K - \Sigma_{KI}\Sigma_{II}^{-1} \vk \mu_I\RT) (\kls{jT} v^{\kls{1}-1/H})} 
}
\COM{=&&
v^{-\abs{I}+1/H-1}
\frac{\mathcal H_{I}(T)}{T} \frac{1}{\sqrt{ (2\pi t_0^{2H})^d \abs{\Sigma}}} e^{-\frac{v^2}{2} g_I(t_0)}\\
&&\times
\int_{-N_{v,M}T v^{1-1/H}}^{(N_{v,M}+1)T v^{1-1/H}}
e^{\frac{-g_I''(t_0) y_v^2}{4}}
\int_{\R^{\abs{I^c}}} e^{-\frac{1}{2 t_0^{2H}} \vk x_{I^c}^ \top (\Sigma^{-1})_{I^cI^c} \vk x_{I^c}}  I_{(\vk x_K <   y_v [(\Sigma^{-1})_{I^cI^c}]^{-1} (\Sigma^{-1} \vk \mu)_{I^c}]_K)} d\vk x_{I^c} dy
\\
\sim
&&
v^{-\abs{I}+1/H-1} \frac{\mathcal H_{I}(T)}{T}\frac{\mcalC_{K,J,M}}{ \sqrt{ (2\pi t_0^{2H} )^d \abs{\Sigma}}} e^{-\frac{v^2}{2} g_I(t_0)}}
\EQNY
Thus, the claim in \eqref{eq:p12v} follows by letting $v\to\IF$ and an application of the Lebesgue dominated convergence theorem. Similarly, we have, as $v\to\IF,$
\BQN\label{eq:p_1v}
p_1(v)\sim
\ksn{ v^{-\abs{I}+1/H-1} \frac{\mathcal H_{I}(T)}{T}\mathcal{C}_K e^{-\frac{v^2}{2} g_I(t_0)}
},  
\EQN
where $\ksn{\mcalC_{K}}$ is given in \eqref{eq:mcalC}.


Next, we consider the term $\pi_M(v)$, where it is sufficient to assume $T$ to be a large number in the sequel. We shall derive a suitable asymptotic upper bound for it, for which we need the following lemma.  Denote
\BQNY
 p(\tau_1,\tau_2;v) = 
\mathbb{P}\Bigg( \exists \begin{array}{l}  t \in [t_0+\tau_1 v^{-1/H},t_0+(\tau_1+1) v^{-1/H}]\\
 s \in [t_0+\tau_2 v^{-1/H},t_0+(\tau_2+1) v^{-1/H}]\end{array}: \vk{X}(t)>(\vk{\nu}+\vk{\mu}t) v,\, \vk{X}(s)>(\vk{\nu}+\vk{\mu}s) v \Bigg).
\EQNY
\BEL
\label{Hsmalldoubleevent}
For any fixed \rJ{$M>0$}, there exist $C_M, v_M>0$  such that, for all $v \geq v_M$,
\BQNY
p(\tau_1,\tau_2;v)\leq C_M \exp\left(
-C_M^{-1} (\tau_2-\tau_1)^{2H}
\right) v^{-|I|} e^{- \frac{v^2 g_I(t_0)}{2}}.
\EQNY
holds uniformly in $\tau_1, \tau_2$ such that $-M v^{1/H-1}\le \tau_1+1\le \tau_2\le M v^{1/H-1}$.
\EEL
\kk{The proof of Lemma \ref{Hsmalldoubleevent} is displayed in Appendix.}

Recall the defination of $p_{j,l;v}$ in \eqref{eq:pjlv}. We shall partition the interval $\Delta_{j;v}$ into segments of the form $[(j+1)T-k-1,(j+1)T-k]$, $0 \leq k \leq \lfloor T\rfloor$ and the interval $\Delta_{l;v}$ into segments of the form $[lT+m,lT+m+1]$, $0\le m \leq \lfloor T\rfloor$. By doing so, we have
\BQNY
p_{j,l;v}\le \sum_{0\le k,m\le \lfloor T\rfloor} p((j+1)T-k-1,lT+m;v).
\EQNY

Applying Lemma \ref{Hsmalldoubleevent} to the above, we have, for all large enough $v$,
\BQNY
p_{j,l;v}
&\leq&
\sum_{
		0 \leq k,m \leq T^{1/3}-1}
C_M \exp\left(
-C_M^{-1} |(l-j-1)T+m+k+1|^{2H}
\right) v^{-|I|} e^{- \frac{v^2 g_I(t_0)}{2}}
\notag
\\
&&+
\underset{T^{1/3}-1 \leq \max(k,m) \leq \lfloor T\rfloor}{\sum_{k,m\ge 0}}
C_M \exp\left(
-C_M^{-1} |(l-j-1)T+m+k+1|^{2H}
\right) v^{-|I|} e^{- \frac{v^2 g_I(t_0)}{2}}
\notag\\
&\le &
T^{2/3}
C_M \exp\left(
-C_M^{-1}( (l-j-1)T)^{2H}
\right) v^{-|I|} e^{- \frac{v^2 g_I(t_0)}{2}}
\notag
\\
&&+ T^2
C_M \exp\left(
-C_M^{-1} ((l-j-1)T+T^{1/3})^{2H}
\right) v^{-|I|} e^{- \frac{v^2 g_I(t_0)}{2}},
\EQNY
\kk{uniformly} for all $j, l$ such that $-N_{v,M} \leq j < l \leq N_{v,M}$. Thus, for all large enough $v,$
\BQN
\label{double_sum_bound}
\pi_M(v)
&\leq&
\sum_{-N_{v,M} \leq j < l \leq N_{v,M}}
p_{j,l;v}\notag
\\
&\leq&
2 N_{v,M}T^{2/3}
C_M \sum_{l \geq 0}  \exp\left(
-C_M^{-1}( (lT)^{2H}
\right) v^{-|I|} e^{- \frac{v^2 g_I(t_0)}{2}}
\\
&&+2 N_{v,M} T^2
C_M \sum_{l \geq 0} \exp\left(
-C_M^{-1} (lT+T^{1/3})^{2H}
\right) v^{-|I|} e^{- \frac{v^2 g_I(t_0)}{2}}.
\notag
\EQN
\kk{Since for all $x \geq 0$ and $T \geq 1$ we have}
\[\exp\left(-C_M^{-1}(lT+x)^{2H}\right)\leq \exp\left(-C_M^{-1}(l+x)^{2H}\right) \leq \rJ{\int_{x+l-1}^{x+l} \exp\left(-C_M^{-1}t^{2H}\right) dt},
\]
\kk{therefore}
\BQNY
\sum_{l \geq 0}
\exp\left(
-C_M^{-1}(lT+x)^{2H}
\right)
\leq \exp\left(
-C_M^{-1}x^{2H}
\right)+
\int_{x}^{\infty} \exp\left(
-C_M^{-1}t^{2H}
\right)dt
\leq \widetilde{C}_M \exp\left(-\widetilde{C}_M^{-1}x^H\right),
\EQNY
\kk{where $\widetilde{C}_M>0$ is a constant independent of $x$.}
Combining this with \eqref{double_sum_bound}, we obtain
\BQN
\label{double_sum_bound_2}
\limsup_{T \to \infty}\rJ{ \lim_{v\to\IF} }\frac{\pi_M(v)}{v^{-|I|+1/H-1}e^{-\frac{v^2 g_I(t_0)}{2}}} = 0.
\EQN
Applying \eqref{eq:p12v}, \eqref{eq:p_1v} and \eqref{double_sum_bound_2} to \eqref{eq:thetaT}, 
we obtain  
\BQN
\ksn{
\mathcal H_{I}\mcalC_{K,M}
}
&\leq&
\liminf_{T \to \infty} \lim_{u\to\IF}\frac{p(v)}{v^{-|I|+1/H-1}e^{-\frac{v^2 g_I(t_0)}{2}}}
\\
&\leq&
\limsup_{T \to \infty} \lim_{u\to\IF} \frac{p(v)}{v^{-|I|+1/H-1}e^{-\frac{v^2 g_I(t_0)}{2}}}
\leq
\ksn{
\mathcal H_{I}\mcalC_{K}
}.
\notag
\EQN
Now, letting $M \to \infty$ in the above, we have
\BQN
\label{final_asymp_1}
p(v) \sim
\ksn{ v^{-\abs{I}+1/H-1}
\mathcal
H_{I}\mcalC_{K} e^{-\frac{v^2}{2} g_I(t_0)}}, \ \ v\to\IF.
\EQN

\underline{Case (ii): \kk{$H>1/2$ and $H \vk{\nu}_{I} \neq (1-H) t_0 \vk{\mu}_{I}$.}} In this case, we shall consider intervals of  length that is shorter than
\kk{the one used in Case (i).}
Precisely, we define, for any fixed integer $T>0$ and $v>0$
\BQNY 
\hDel_{j;v}=\hDel_{j;v}(T)= [t_0+j Tv^{-2}, t_0+(j+1) Tv^{-2}],\ \ \ -\hN_v-1\le j\le \hN_v,
\EQNY
where $\hN_v=\lfloor T^{-1} \ln(v) v\rfloor$. Also denote $\hN_{v,M} = \lfloor T^{-1} M v \rfloor$.
By Bonferroni's inequality we have
\BQN\label{eq:thetaT2}
\hp_1(v)\ge p(v) \ge \hp_{2,M}(v)-\hpi_M(v),
\EQN
where
\begin{equation*} 
\hp_1(v)= \sum_{j=-\hN_v-1}^{\hN_v}\hp_{j;v},\ \
\hp_{2,M}(v) = \sum_{j=-\hN_{v,M}}^{\hN_{v,M}}\hp_{j;v},\ \
\hpi_M(v) = \sum_{ -\hN_{v,M}\le j< l\le  \hN_{v,M}}\hp_{j,l;v},
\end{equation*}
with
\begin{equation*} 
\hp_{j;v}=\pk{\exists_{t\in\hDel_{j;v}} \vk{X}(t)>
(\vk{\alpha}+\vk{\mu} t)v}
\end{equation*}
and
$$
\hp_{j,l;v}=\pk{\exists_{t\in \hDel_{j;v}}  \vk{X}(t)> (\vk{\alpha}+\vk{\mu} t)v,
 \ \exists_{t\in \hDel_{l;v}} \vk{X}(t)> (\vk{\alpha}+\vk{\mu} t)  v  }.
 $$

Similarly as in Case (i), we shall deal with the single-sum $\hp_1(v)$, $\hp_{2,M}(v)$ and the double-sum $\hpi_M(v)$, respectively. For the asymptotics of the single-sum terms, we shall use the following uniform version of a generalized Pickands lemma evaluated on a shorter interval. The proof of Lemma \ref{lem:Pick2} is deferred to Appendix.

\BEL \label{lem:Pick2} Fix $T>0$.  We have, as $v\to\IF,$
\BQNY
&&\pk{\exists_{t \in[t_0+\tau v^{-2}, t_0+(\tau+T) v^{-2}]}
 \vk X(t) > (\vk \nu +\vk \mu t) v }\\
&&\sim  v^{-\abs{I}} \frac{ t_0^{2H\abs{I}}}{  \prod_{i \in I}w_i } \LT(1+\frac{T}{t_0^{2H}}\sum_{i \in I} \left(w_i \mu_i - \frac{H}{t_0} w_i b_i\right)_{-} \RT) e^{-\frac{v^2}{2} g_I(t_0+\tau v^{-2})}\\
&&
\quad
\ksn{
	\times \pk{\vk Y_K<- t_0^{-H} \LT(\vk \mu_K - \Sigma_{KI}\Sigma_{II}^{-1} \vk \mu_I\RT) (\tau v^{-1})}  /\sqrt{ (2\pi  t_0^{2H} )^{\abs{I}} }/\sqrt{ \abs{\Sigma_{II}}}
},
\EQNY
holds uniformly in $\tau$ such that $\abs{\tau}<  T ( \hN_v+1)$.
\EEL

With Lemma \ref{lem:Pick2} given above, it follows \kk{by}  the same lines of reasoning as in Case (i) that
\BQN
\begin{aligned}\label{eq:p14v}
\hp_{2,M}(v)&\sim
v^{-\abs{I}+1}
\ksn{
\frac{t_0^{2H|I|} \mcalC_{K,M}}{T \prod_{i \in I} w_i}
} \LT(1+\frac{T}{t_0^{2H}}\sum_{i \in I} \left(w_i \mu_i - \frac{H}{t_0} w_i b_i\right)_{-} \RT) e^{-\frac{v^2}{2} g_I(t_0)},
\\
\hp_1(v) &\sim  v^{-\abs{I}+1}
\ksn{
\frac{t_0^{2H|I|} \mcalC_{K}}{T \prod_{i \in I} w_i}
} \LT(1+\frac{T}{t_0^{2H}}\sum_{i \in I} \left(w_i \mu_i - \frac{H}{t_0} w_i b_i\right)_{-} \RT) e^{-\frac{v^2}{2} g_I(t_0)}
\end{aligned}
\EQN
hold, as $v\to\IF$. 


Next, in order to derive a suitable upper bound for the double-sum term $\hpi_M(u)$, we need an analogue of Lemma \ref{Hsmalldoubleevent}.  It is worth noting that Lemma \ref{Hbigdoubleevent} below looks similar to Lemma \ref{Hsmalldoubleevent}, but the approach used to prove it is quite different, which
is displayed in Appendix.

Denote
\BQNY
 \hp(\tau_1,\tau_2;v) = 
\mathbb{P}\Bigg( \exists \begin{array}{l}  t \in [t_0+\tau_1 v^{-2},t_0+(\tau_1+1) v^{-2}]\\
 s \in [t_0+\tau_2 v^{-2},t_0+(\tau_2+1) v^{-2}]\end{array}: \vk{X}(t)>(\vk{\nu}+\vk{\mu}t) v,\, \vk{X}(s)>(\vk{\nu}+\vk{\mu}s) v \Bigg).
\EQNY
\BEL
\label{Hbigdoubleevent}
For any fixed \rJ{$M>0$}, there exist $C_M, v_M>0$  such that, for all $v \geq v_M$,
\BQNY
\hp(\tau_1,\tau_2;v)\leq C_M \exp\left(
-C_M^{-1} (\tau_2-\tau_1)
\right) v^{-|I|} e^{- \frac{v^2 g_I(t_0)}{2}}
\EQNY
holds uniformly in $\tau_1, \tau_2$ such that $-M v\le \tau_1+1\le \tau_2\le M v$.
\EEL

Similarly to \eqref{double_sum_bound_2}, we have from the above lemma that
\BQN\label{eq:hpi_M}
\limsup_{T \to \infty}\rJ{ \lim_{v\to\IF} }\frac{\hpi_M(v)}{v^{-|I|+1}e^{-\frac{v^2 g_I(t_0)}{2}}} = 0.
\EQN
Consequently,
we conclude from \eqref{eq:thetaT2}-\eqref{eq:hpi_M} that 
\BQN
\label{final_asymp_2}
p(v) \sim
v^{-\abs{I}+1}
\ksn{ \frac{t_0^{2H(|I|-1)} \mcalC_{K}}{ \prod_{i \in I} w_i}
} \LT(\sum_{i \in I} \left(w_i \mu_i - \frac{H}{t_0} w_i b_i\right)_{-} \RT) e^{-\frac{v^2}{2} g_I(t_0)}, \ \ v\to\IF.
\EQN


\subsection*{\underline{Step 2: Analysis of $\Pi(v)$}} \kk{In order to obtain a sharp upper bound for $\Pi(v)$,
we can adapt the same arguments as in  Lemma 4.1 of \cite{DHJRBM}, which gives that, for sufficiently large $v$ it holds that}
\BQN \label{eq:Piv}
\begin{aligned}
\Pi(v)&=\pk{\exists_{t \in\widetilde \Delta_v}  \vk X(t) > (\vk \nu +\vk \mu t) v  }\\
&\le  C_0 \exp\LT(- \frac{v^2}{2} g_I(t_0)- \frac{\min(g''(t_0+), g''(t_0-))}{C_1} (\ln(v))^2\RT),
\end{aligned}
\EQN
where \kls{$C_0, C_1>0$ are some constants} and $g''(t_0\pm)$ are the one-side second derivatives of $g(t)$, with $g''(t_0\pm)>0$ that can be confirmed  as in the proof of Lemma \ref{lem:gc1}. It is noted that generally $g''(t_0+)\neq g''(t_0-) $; see Remark 5.7 of \cite{DHJRBM}   for such an example.

\medskip

\kk{In order to complete the proof of Theorem \ref{Thm:main}, we note that
by combining  \eqref{eq:Piv} with \eqref{final_asymp_1} or \eqref{final_asymp_2}
we get
that
\[ \Pi(v)=o(p(v)),
\]
as $v\to\infty$. The above, together with
\eqref{eq:pPpv}  (recalling that $v=u^{1-H}$) completes the proof.} \QED

\section{Appendix: Additional proofs}

{\bf Proof of Lemma \ref{lem:gc1}.} \kk{The proof follows by the use of similar arguments as in the proof of} Lemma 2.2 in \cite{DHJRBM}. For completeness we present
\kk{the main steps of argumentation} and only highlight the key differences. First, 
 note that $h(t)=  \inf_{\vk{v} \ge \vk\nu+\vk{\mu} t}  \vk{v}^\top \Sigma^{-1}  \vk{v} \in C^1(0,\IF)$ has been proved in \cite{DHJRBM}, thus $g\in C^1(0,\IF)$ is established. Next,
 denote $I(t)\subseteq \{1,2,\ldots, d\}$ to be the index set of the solution to the quadratic programming problem $P_{\Sigma}(\vk b(t))$ for any fixed $t>0$. It follows from Lemma A.4 of  \cite{DHJRBM} that $I(t), t>0$ is an almost piecewise constant set function. Namely,
$$I(t)=\sum_{j} I_j\mathbb{I}_{(t\in U_j)},$$
where 
 $U_j$'s are of the following form
 \begin{equation*} 
 (a,b), [a,b), (a,b], [a,b], \{a\}, (b,\IF), [b,\IF),
 \end{equation*}
where $0<a<b< \infty$
and $I_j\subseteq \{1,\ldots,d\}$. Therefore,
\BQNY
g(t)=g_{I_j}(t)=\frac{\vk\mu_{I_j}^\top \Sigma^{-1}_{I_jI_j}  \vk\mu_{I_j}  t^2 +2\vk \nu_{I_j}^\top \Sigma^{-1}_{I_jI_j}  \vk\mu_{I_j}  t +\vk\nu_{I_j}^\top \Sigma^{-1}_{I_jI_j}  \vk\nu_{I_j}  }{t^{2H}}, \quad t \in U_j^o,
\EQNY
where $U_j^o$ is the inner set of $U_j$.
Furthermore, for any fixed $I_j$, the first derivative of $g_{I_j}$, 
\BQNY
g_{I_j}'(t)= \frac{2 \vk\mu_{I_j}^\top \Sigma^{-1}_{I_jI_j}  \vk\mu_{I_j} (1-H) t^2 +2(1-2H)\vk \nu_{I_j}^\top \Sigma^{-1}_{I_jI_j}  \vk\mu_{I_j}  t -2H \vk\nu_{I_j}^\top \Sigma^{-1}_{I_jI_j}  \vk\nu_{I_j}  }{t^{2H+1}}. 
\EQNY
is negative on the left of the positive root of $g_{I_j}'(t)=0$ and then becomes positive on the right of this root. This means that the function $g_{I_j}(t), t>0$ is decreasing to the left of some point and then becomes increasing. Next, we have $g(t)\to\IF$ as $t\to\IF$ or $t\to 0$. From these and the fact that $g\in C^1(0,\IF)$, and using the same arguments as in Lemma 2.2 of \cite{DHJRBM} we can conclude that  the minimizer of the function $g(t), t>0$ is given by  $t_0 \in U_j$ for some $j$, which \ccP{must} be of the form \eqref{eq:t0} and satisfies \eqref{eq:intr1}. Moreover, elementary calculations show that \eqref{eq:gt0pm} is valid, where $g''_I(t_0)>0$ follows from the fact that $2H g_I'(t)+t g_I''(t)>0$ for any non-empty set $I \subset \{1,\dots,d\}$ and
$t>0$. The proof is complete. \QED

{\bf Proof of Lemma \ref{lem:Pick1}.}   It follows that
\BQNY
&&\pk{\exists_{t \in[t_0+\tau v^{-1/H}, t_0+(\tau+T) v^{-1/H}]}
 \vk X(t) > (\vk \nu +\vk \mu t) v }\\
&&= \pk{\exists_{t \in[0, T]}
 \vk X(t_0+\tau v^{-1/H}+tv^{-1/H}) > (\vk \nu +\vk \mu (t_0+\tau v^{-1/H}+t v^{-1/H})) v }\\
&&=  \pk{\exists_{t \in[0, T]}
\vk X_{v,\tau}(t) > \vk b v + \tau \vk \mu v^{1-1/H}+t \vk \mu v^{1-1/H} },
\EQNY
where $\vk X_{v,\tau}(t) = \vk X(t_0+\tau v^{-1/H}+t v^{-1/H}).$ We shall follow some ideas in the proof of Lemma 4.7 in \cite{DHW}.
Let $\widetilde{\vk b}$ be the optimal solution of the optimization problem $P_\Sigma (\vk b)$. 
 We have
\BQNY
&& \pk{\exists_{t \in[0, T]}
\vk X_{v,\tau}(t) > \vk b v + \tau \vk \mu v^{1-1/H}+t \vk \mu v^{1-1/H} }\\
&&= \pk{\exists_{t \in[0, T]}
\vk X_{v,\tau}(t) - \wtb v> (\vk b- \wtb) v + \tau \vk \mu v^{1-1/H}+t \vk \mu v^{1-1/H} }.
\EQNY
Define
$$
\vk Z_{v,\tau}(t) :=\overline{\vk v}(\vk X_{v,\tau}(t) - \wtb v- \tau \vk \mu v^{1-1/H}-t \vk \mu v^{1-1/H})+\vk x,
$$
where $\olv$ has all components equal to $v$ for the indices in $I$, and 1 for the indices in $I^c=K\cup J$.
It then follows that
\BQNY
&&  \pk{\exists_{t \in[0, T]}
\vk X_{v,\tau}(t) - \wtb v> (\vk b- \wtb) v + \tau \vk \mu v^{1-1/H}+t \vk \mu v^{1-1/H} }\\
&& = v^{-\abs{I}} \int_{\R^d}  \pk{\exists_{t \in[0, T]}
\vk Z_{v,\tau}(t) > (\vk b- \wtb) \olv v+\vk x \mid  \vk Z_{v,\tau}(0)=\vk 0} \varphi_{\Sigma_{v,\tau}}(\wtb  v+ \tau \vk \mu v^{1-1/H}-\vk x/\olv) d\vk x,
\EQNY
with
$$
\Sigma_{v,\tau} := \E{\vk X_{v,\tau}(0) \vk X_{v,\tau}(0)^\top} = (t_0+\tau v^{-1/H})^{2H} \Sigma.
$$
Further, denote
$$
\vk \chi_{v,\tau}(t) := (\vk Z_{v,\tau}(t) \mid \vk Z_{v,\tau}(0)=\vk 0)
$$
and
$$
r_{v,\tau}(t,s):=\frac{1}{2}\LT( (t_0+\tau v^{-1/H}+tv^{-1/H})^{2H} + (t_0+\tau v^{-1/H}+sv^{-1/H})^{2H} -\abs{t-s}^{2H} v^{-2}\RT).
$$
We derive that  
\BQNY
\E{\vk \chi_{v,\tau}(t) }&=& \overline{\vk v}\LT( \E{\vk X_{v,\tau}(t)\mid \vk X_{v,\tau}(0)= \wtb v+ \tau \vk \mu v^{1-1/H}-\vk x/\olv } - \wtb v- \tau \vk \mu v^{1-1/H}-t \vk \mu v^{1-1/H}\RT)+\vk x \\
&=& \overline{\vk v}\LT( r_{v,\tau}(t,0) r_{v,\tau}(0,0)^{-1}( \wtb v+ \tau \vk \mu v^{1-1/H}-\vk x/\olv)  - \wtb v- \tau \vk \mu v^{1-1/H}-t \vk \mu v^{1-1/H}\RT)+\vk x \\
&=& \LT(1-  r_{v,\tau}(t,0) r_{v,\tau}(0,0)^{-1}\RT) \vk x + \overline{\vk v}\LT( ( r_{v,\tau}(t,0) r_{v,\tau}(0,0)^{-1} -1) ( \wtb v+ \tau \vk \mu v^{1-1/H} )  -t \vk \mu v^{1-1/H}\RT).
\EQNY
Next, it follows that
\BQNY
 r_{v,\tau}(t,0) r_{v,\tau}(0,0)^{-1} -1 = \frac{(t_0+\tau v^{-1/H}+tv^{-1/H})^{2H} - (t_0+\tau v^{-1/H})^{2H}}{2 (t_0+\tau v^{-1/H})^{2H}} -  \frac{t^{2H} v^{-2}}{2 (t_0+\tau v^{-1/H})^{2H}},
\EQNY
and
\BQNY
(t_0+\tau v^{-1/H}+tv^{-1/H})^{2H} - (t_0+\tau v^{-1/H})^{2H} = 2H (t_0+\tau v^{-1/H})^{2H-1} tv^{-1/H} + H(2H-1) (t_0+\tau' v^{-1/H})^{2H-2} (tv^{-1/H})^2, 
\EQNY
with some $\tau' \in [\tau, \tau+t]$. Some elementary calculations yield that, as $v\to\IF$,
\BQN
\label{eq:chi_t}
\E{\vk \chi_{v,\tau}(t) }
\to  \begin{pmatrix}
		-\frac{1}{2t_0^{2H}}\vk b_I t^{2H} \\
		\vk 0 _{I^c}
	\end{pmatrix}
\EQN
holds uniformly for  $\tau$ such that $\abs{\tau}\leq T(N_v+1)$, where when $H>1/2$, the condition $H \vk \nu_I = (1-H) t_0 \vk \mu_I$ was used.

Now, we analyze the covariance function of $\vk \chi_{v,\tau}(t), t\ge 0$
\BQNY
\E{[\vk \chi_{v,\tau}(t) -\E{\vk \chi_{v,\tau}(t)}][\vk \chi_{v,\tau}(s) -\E{\vk \chi_{v,\tau}(s)}]^\top}
=\text{diag} (\olv)          \LT(  r_{v,\tau}(t,s) -\frac{r_{v,\tau}(t,0)  r_{v,\tau}(0,s)  }{ r_{v,\tau}(0,0)}   \RT)     \Sigma  \  \text{diag} (\olv).
\EQNY
Note that
\BQNY
&& r_{v,\tau}(t,s) r_{v,\tau}(0,0)  -r_{v,\tau}(t,0)  r_{v,\tau}(0,s) \\
&&= ( r_{v,\tau}(t,s) -  r_{v,\tau}(t,0) ) r_{v,\tau}(0,0) + r_{v,\tau}(t,0) (  r_{v,\tau}(0,0)- r_{v,\tau}(0,s) )\\
&&=\frac{1}{2}\LT( (t_0+\tau v^{-1/H}+sv^{-1/H})^{2H} - (t_0+\tau v^{-1/H})^{2H} \RT)( r_{v,\tau}(0,0)- r_{v,\tau}(t,0) )  \\
&& \ \ \ \ + \frac{1}{2} \LT(t^{2H} r_{v,\tau}(0,0) +s^{2H}  r_{v,\tau}(t,0) - \abs{t-s}^{2H}  r_{v,\tau}(0,0) \RT) v^{-2}.
\EQNY
Thus, similarly to the calculations for the mean, we obtain that, as $v\to\IF,$
\BQNY
\E{[\vk \chi_{v,\tau}(t) -\E{\vk \chi_{v,\tau}(t)}][\vk \chi_{v,\tau}(s) -\E{\vk \chi_{v,\tau}(s)}]^\top}\to  \begin{pmatrix}
	K(t,s)\Sigma_{II} &\vk 0_{II^c}\\
		\vk 0_{I^c I} &  \vk 0_{I^c I^c}
	\end{pmatrix}
\EQNY
holds uniformly for   $\tau$ such that $\abs{\tau} \leq T(N_v+1)$, where $K(t,s)=\frac{1}{2}(t^{2H}+s^{2H}-\abs{t-s}^{2H})$.
Next, we show that the process $\{\vk{\chi}_{v,\tau}(t), t\in[0,T]\}$ is tight for all $\tau$ such that $|\tau| \leq T(N_v+1)$ and large enough $v$. To this end, it is sufficient to show the tightness of the conditional process $\{\overline{\vk v}\vk X_{v,\tau}(t) \mid \vk Z_{v,\tau}(0)=\vk 0, t\in [0,T]\}$.
Let us note that for any  jointly Gaussian distributed  random vectors $\vk U$ and  $\vk Y$, it is known that
\BQN
\label{cond}
\VA{\vk U | \vk Y =\vk y} \leq \VA{\vk U},
\EQN
where
$\VA{\vk V} := \E{(\vk V-\E{\vk V})^{\top}(\vk V-\E{\vk V})}$
\kls{is the sum of variances of the random (column) vector $\vk V$}.
Hence, for all $t,s\in[0,T],$
\BQNY
\VA{\overline{\vk v}\vk X_{v,\tau}(t) - \overline{\vk v}\vk X_{v,\tau}(s)   \mid \vk Z_{v,\tau}(0)=\vk 0}\le \VA{\overline{\vk v}\vk X_{v,\tau}(t) - \overline{\vk v}\vk X_{v,\tau}(s)}
\le C |t-s|^{2H}
\EQNY
holds uniformly for all $\tau$ such that $|\tau| \leq T(N_v+1)$ and large enough $v$, where $C>0$ is a constant.
Thus,  $\{\vk{\chi}_{v,\tau}(t), t\in[0,T]\}$ is tight. 

Therefore, following  similar arguments as in \cite{DHW} we conclude that, as $v\to\IF$,
\BQNY
  \pk{\exists_{t \in[0, T]}
\vk Z_{v,\tau}(t) > (\vk b- \wtb) \olv v +\vk x \mid  \vk Z_{v,\tau}(0)=\vk 0}  \to    \pk{\exists_{t \in[0, T]}
\vk W_{I}(t)  -  \frac{1}{2t_0^{2H}}\vk b_I t^{2H}>\vk x_I }\cdot  \kls{\mathbb {I}}_{(\vk x_K < \vk 0_K)}
\EQNY
uniformly for all $\tau$ such that $|\tau| \leq T(N_v+1)$, where $\vk W_{I}(t) =D \vk B_{H,I}(t)$ with $D D^\top = \Sigma_{II}$.

Next, we have
\BQNY
 &&\varphi_{\Sigma_{v,\tau}}(\wtb  v + \tau \vk \mu v^{1-1/H}-\vk x/\olv) =\frac{1}{\sqrt{(2\pi)^d \abs{\Sigma_{v,\tau}}} }\\
&& \ \ \ \ \ \times \exp\LT(- \frac{1}{2 (t_0+\tau v^{-1/H})^{2H}} (\wtb  v + \tau \vk \mu v^{1-1/H}-\vk x/\olv) ^\top\Sigma^{-1} (\wtb  v + \tau \vk \mu v^{1-1/H}-\vk x/\olv) \RT)\\
&&=  \varphi_{\Sigma_{v,\tau}}(\wtb  v + \tau \vk \mu v^{1-1/H}) \exp\LT(\frac{1}{ (t_0+\tau v^{-1/H})^{2H}} (\vk x/\olv)^\top \Sigma^{-1}  (\wtb  v + \tau \vk \mu v^{1-1/H}) -  \frac{1}{2 (t_0+\tau v^{-1/H})^{2H}}  (\vk x/\olv)^\top \Sigma^{-1}  (\vk x/\olv) \RT).
\EQNY
\kk{Since
$
\vk w = \Sigma^{-1}\wtb
$
and (by \eqref{eq:ww})  $\vk w_I =(\Sigma_{II})^{-1} \vk b_I>0$ and $\vk w_{I^c} =\vk 0_{I^c}$,}
the above exponent is asymptotically equal to, as $v\to\IF,$
\BQNY
\frac{1}{t_0^{2H}} \vk w_I^\top \vk x_I +\frac{1}{t_0^{2H}} \vk x_{I^c}^\top (\Sigma^{-1}  \vk \mu)_{I^c} ( \tau v^{1-1/H})-\frac{1}{2 t_0^{2H}} \vk x_{I^c}^ \top (\Sigma^{-1})_{I^cI^c} \vk x_{I^c}
\EQNY
uniformly for  $\tau$ such that $\abs{\tau}\leq T(N_v+1)$.
Next, we shall rewrite $ \varphi_{\Sigma_{v,\tau}}(\wtb  v + \tau \vk \mu v^{1-1/H})$.  
 Note that for any $\vk y\in \R^d$, we have
\BQNY
(\wtb +\vk y)^\top \Sigma^{-1} (\wtb +\vk y) &=& \vk b_I^\top (\Sigma_{II})^{-1} \vk b_I + 2 \vk b_I^\top (\Sigma_{II})^{-1} \vk y_I +\vk y^\top \Sigma^{-1}  \vk y\\
&=& (\wtb +\vk y)_I^\top (\Sigma_{II})^{-1} (\wtb +\vk y)_I  -  \vk y_I^\top (\Sigma_{II})^{-1}  \vk y_I +\vk y^\top \Sigma^{-1}  \vk y.
\EQNY
Therefore, as $v\to\IF$,
\BQNY
\varphi_{\Sigma_{v,\tau}}(\wtb  v + \tau \vk \mu v^{1-1/H})& \sim &  \frac{1}{ \sqrt{ (2\pi )^d \abs{\Sigma_{v,\tau}}}}\exp\LT(-\frac{v^2}{2} g_I(t_0 +\tau v^{-1/H})\RT) \\
&&\times \exp\LT(- \frac{(\tau v^{1-1/H})^2}{2t_0^{2H}} \LT( \vk \mu^\top \Sigma^{-1}  \vk \mu-  \vk \mu_I^\top (\Sigma_{II})^{-1}  \vk \mu_I  \RT)\RT)
\EQNY
uniformly for  $\tau$ such that $\abs{\tau}\leq T(N_v+1)$.
Putting everything together and using the dominated convergence theorem as in \cite{DHW}  (\kk{we omit the standard details}), we obtain, as $v\to\IF$,
\BQN\label{eq:Pick1}
&&\pk{\exists_{t \in[t_0+\tau v^{-1/H}, t_0+(\tau+T) v^{-1/H}]}
 \vk X(t) > (\vk \nu +\vk \mu t) v }\nonumber \\
&&\sim  v^{-\abs{I}} \frac{\mathcal H_{I}(T)}{\sqrt{ (2\pi t_0^{2H})^d \abs{\Sigma}}} e^{-\frac{v^2}{2} g_I(t_0+\tau v^{-1/H})} \times \exp\LT(- \frac{(\tau v^{1-1/H})^2}{2t_0^{2H}} \LT( \vk \mu^\top \Sigma^{-1}  \vk \mu-  \vk \mu_I^\top (\Sigma_{II})^{-1}  \vk \mu_I  \RT)\RT) \nonumber\\
&&\ \ \times \int_{\R^{\abs{I^c}}} e^{\frac{1}{t_0^{2H}} \vk x_{I^c}^\top (\Sigma^{-1}  \vk \mu)_{I^c} ( \tau v^{1-1/H})-\frac{1}{2 t_0^{2H}} \vk x_{I^c}^ \top (\Sigma^{-1})_{I^cI^c} \vk x_{I^c}}  \kls{\mathbb {I}}_{(\vk x_K < \vk 0_K)} d\vk x_{I^c}
\EQN
uniformly for  $\tau$ such that $\abs{\tau}\leq T(N_v+1)$.
Furthermore, using the Schur complement of invertible block matrix and some elementary calculations, it can be derived that
\BQNY
\vk \mu^\top \Sigma^{-1}  \vk \mu-  \vk \mu_I^\top (\Sigma_{II})^{-1}  \vk \mu_I   = (\Sigma^{-1} \vk \mu)_{I^c}^\top[(\Sigma^{-1})_{I^cI^c}]^{-1} (\Sigma^{-1} \vk \mu)_{I^c}.
\EQNY
Therefore,
\BQNY
&&- \frac{(\tau v^{1-1/H})^2}{2t_0^{2H}} \LT( \vk \mu^\top \Sigma^{-1}  \vk \mu-  \vk \mu_I^\top (\Sigma_{II})^{-1}  \vk \mu_I  \RT)+\frac{1}{t_0^{2H}} \vk x_{I^c}^\top (\Sigma^{-1}  \vk \mu)_{I^c} ( \tau v^{1-1/H})-\frac{1}{2 t_0^{2H}} \vk x_{I^c}^ \top (\Sigma^{-1})_{I^cI^c} \vk x_{I^c}\\
&&=-\frac{1}{2t_0^{2H}}[\vk x_{I^c} -   ( \tau v^{1-1/H}) [(\Sigma^{-1})_{I^cI^c}]^{-1} (\Sigma^{-1} \vk \mu)_{I^c}]^\top (\Sigma^{-1})_{I^cI^c} [\vk x_{I^c} -   ( \tau v^{1-1/H}) [(\Sigma^{-1})_{I^cI^c}]^{-1} (\Sigma^{-1} \vk \mu)_{I^c}].
\EQNY
\kk{
Hence
\begin{eqnarray*}
\lefteqn{
\exp\LT(- \frac{(\tau v^{1-1/H})^2}{2t_0^{2H}} \LT( \vk \mu^\top \Sigma^{-1}  \vk \mu-  \vk \mu_I^\top (\Sigma_{II})^{-1}  \vk \mu_I  \RT)\RT)
\times \int_{\R^{\abs{I^c}}} e^{\frac{1}{t_0^{2H}} \vk x_{I^c}^\top (\Sigma^{-1}  \vk \mu)_{I^c} ( \tau v^{1-1/H})-\frac{1}{2 t_0^{2H}} \vk x_{I^c}^ \top (\Sigma^{-1})_{I^cI^c} \vk x_{I^c}}   \kls{\mathbb {I}}_{(\vk x_K < \vk 0_K)} d\vk x_{I^c}}\\
&=&
\int_{\R^{\abs{I^c}}}
 e^{-\frac{1}{2t_0^{2H}}[\vk x_{I^c} -   ( \tau v^{1-1/H}) [(\Sigma^{-1})_{I^cI^c}]^{-1} (\Sigma^{-1} \vk \mu)_{I^c}]^\top (\Sigma^{-1})_{I^cI^c} [\vk x_{I^c} -   ( \tau v^{1-1/H}) [(\Sigma^{-1})_{I^cI^c}]^{-1} (\Sigma^{-1} \vk \mu)_{I^c}]}
  \kls{\mathbb {I}}_{(\vk x_K < \vk 0_K )} d\vk x_{I^c}\\
&=&
\int_{\R^{\abs{I^c}}} e^{-\frac{1}{2 t_0^{2H}} \vk x_{I^c}^ \top (\Sigma^{-1})_{I^cI^c} \vk x_{I^c}}   \kls{\mathbb {I}}_{(\vk x_K <  - ( \tau v^{1-1/H}) [[(\Sigma^{-1})_{I^cI^c}]^{-1} (\Sigma^{-1} \vk \mu)_{I^c}]_K)} d\vk x_{I^c}.
\end{eqnarray*}
}
\kls{On the other hand, we have, from the Schur complement formula,
\BQNY
&&[(\Sigma^{-1})_{I^cI^c}]^{-1}= \Sigma_{I^c I^c} -  \Sigma_{I^c I}\Sigma_{II}^{-1} \Sigma_{I I^c}, \\
&&(\Sigma^{-1})_{I^cI}=- (\Sigma^{-1})_{I^cI^c}  \Sigma_{I^c I}\Sigma_{II}^{-1}.
\EQNY
}
This implies that
\ksn{\BQN
\label{integr_expr}
&&
\int_{\R^{\abs{I^c}}} e^{-\frac{1}{2 t_0^{2H}} \vk x_{I^c}^ \top (\Sigma^{-1})_{I^cI^c} \vk x_{I^c}}   \kls{\mathbb {I}}_{(\vk x_K <  - ( \tau v^{1-1/H}) [[(\Sigma^{-1})_{I^cI^c}]^{-1} (\Sigma^{-1} \vk \mu)_{I^c}]_K)} d\vk x_{I^c}
\notag
\\
&&
=
\sqrt{\frac{(2\pi t_0^{2H})^{|I^c|}}
{
|(\Sigma^{-1})_{I^c I^c}|
}
}
\  \pk{\vk Z_{K}<- ( \tau v^{1-1/H}) t_0^{-H} \LT(\vk \mu_{I^c} - \Sigma_{I^c I}\Sigma_{II}^{-1} \vk \mu_I\RT)_K, \vk Z_{I^c \setminus K}<\vk\IF_{I^c \setminus K}}
\notag
\\
&&
= \sqrt{\frac{(2\pi t_0^{2H})^{|I^c|}}
{
|(\Sigma^{-1})_{I^c I^c}|
}
} \ \pk{\vk Y_K< - ( \tau v^{1-1/H}) t_0^{-H} \LT(\vk \mu_{K} - \Sigma_{K I}\Sigma_{II}^{-1} \vk \mu_I \RT) }, 
\EQN
where 
\[
\vk{Z}_{I^c} \overset{d}\sim \mathcal N\LT(\vk 0_{I^c}, \Sigma_{I^c I^c} -  \Sigma_{I^c I}\Sigma_{II}^{-1} \Sigma_{I I^c} \RT),\ \ \vk{Z}_K \overset{d}\sim
\vk{Y}_K \overset{d}\sim \mathcal N\LT(\vk 0_K, \Sigma_{K K} -
\Sigma_{K I}\Sigma_{II}^{-1} \Sigma_{I K}
 \RT).
\]
\COM{ and the facts that
$\LT(\vk \mu_{I^c} - \Sigma_{I^c I}\Sigma_{II}^{-1} \vk \mu_I\RT)_K=\vk \mu_{K} - \Sigma_{K I}\Sigma_{II}^{-1} \vk \mu_I$ and
\[
\vk{Z}_K \overset{d}\sim \mathcal N\LT(\vk 0_K, \mathrm{Id}_{K I^c}\Sigma_{I^c I^c} \mathrm{Id}_{I^c K} -
\mathrm{Id}_{K I^c} \Sigma_{I^c I}\Sigma_{II}^{-1} \Sigma_{I I^c}
\mathrm{Id}_{I^c K} \RT)
=
\mathcal N\LT(\vk 0_K, \Sigma_{K K} -
\Sigma_{K I}\Sigma_{II}^{-1} \Sigma_{I K}
 \RT),
\]
where for $I,J \subset
\{1,\dots,d\}$ by $\mathrm{Id}_{IJ}$ we denote the sub-block of the $d \times d$ identity matrix determined by $I$ and $J$.}
Consequently,  the required result is obtained by noting that  $|\Sigma||(\Sigma^{-1})_{I^c I^c}|=|\Sigma_{II}|$ (which follows from the Schur complement formula). \QED
}

{\bf Proof of Lemma \ref{Hsmalldoubleevent}:}
Hereafter the exact values of the constants  $C_M,v_M,C_{M,\epsilon}$ and $v_{M,\epsilon}$ are not important and can be changed \rJ{from} line to line, where \rJ{$\epsilon>0$ is a small enough constant which may also \kk{change} from line to line}. Moreover,  all the inequalities hold uniformly in $\tau_1,\tau_2 \in [-Mv^{1/H-1}-1,Mv^{1/H-1}]$ such that \rJ{$\tau_2 -\tau_1 \geq 1$ for large enough $v$}, and the constants do not depend on $\tau_1,\tau_2$.

Without loss of generality, we shall assume that $I = \{1,\ldots, d\},$
since otherwise, an upper bound for $p(\tau_1, \tau_2; v)$ with only $I$ components can be used.  For convenience of notation 
we sometimes also keep index $I$ though it can be omitted or replaced by $ \{1,\ldots, d\}$.
Denote
\BQNY
\widetilde{\vk{X}}_v(t,s):=
\frac{1}{2}(\vk{X}_v(t)+\vk{X}_v(s)), \ \ \ \text{with}\  \vk{X}_v(t):=
\frac{\vk X(t_0+tv^{-1/H})}{(t_0+tv^{-1/H})^H}
\EQNY
and
\BQNY
\vk{b}_v(t_1,t_2) :=
\frac{1}{2}
\left(
\vk{\zeta}(t_0+t_1v^{-1/H})
+
\vk{\zeta}(t_0+t_2v^{-1/H})
\right), \ \ \ \text{with}\ \vk{\zeta}(t):= \frac{\vk{\mu}+\vk{\nu}t}{t^H}.
\EQNY
It follows that
\BQN\label{eq:p_tau_v}
p(\tau_1,\tau_2;v) \leq \pk{\exists t_1 \in [\tau_1,\tau_1+1],t_2 \in [\tau_2,\tau_2+1]:\widetilde{\vk{X}}_v(t_1,t_2)>v\vk{b}_v(t_1,t_2)}.
\EQN
Next, \kk{let} 
\BQNY
r_v(t,s) :=
\E{
	\frac{B_{H,1}(t_0+tv^{-1/H})}{(t_0+tv^{-1/H})^H}
	\frac{B_{H,1}(t_0+sv^{-1/H})}{(t_0+sv^{-1/H})^H}
}
=
\frac{(t_0+tv^{-1/H})^{2H}+(t_0+sv^{-1/H})^{2H}-|t-s|^{2H} v^{-2}}{2(t_0+tv^{-1/H})^{H}(t_0+sv^{-1/H})^{H}}.
\EQNY
We have, for $t_1,s_1 \in [\tau_1,\tau_1+1],t_2,s_2 \in [\tau_2,\tau_2+1]$,
\BQNY
\widetilde{R}_v(t_1,t_2,s_1,s_2)& :=& \E{\widetilde{\vk{X}}_v(t_1,t_2)\widetilde{\vk{X}}^{\top}_v(s_1,s_2)}\\
&=&\frac{1}{4}  \left(R_v(t_1,s_1)+R_v(t_1,s_2)+R_v(t_2,s_1)+R_v(t_2,s_2)\right)\\
&=& \widetilde{r}_v(t_1,t_2,s_1,s_2)  \Sigma,
\EQNY
where
\BQNY
&& R_v(t,s) := \E{\vk{X}_v(t)\vk{X}^{\top}_v(s)} =r_v(t,s)\Sigma,\\
&&  \widetilde{r}_v(t_1,t_2,s_1,s_2):=\frac{1}{4}  \left(r_v(t_1,s_1)+r_v(t_1,s_2)+r_v(t_2,s_1)+r_v(t_2,s_2)\right).
\EQNY
Moreover, we denote
\BQNY
\widetilde{\Sigma}_v(t_1,t_2) := \widetilde{R}(t_1,t_2,t_1,t_2), \ \ \widetilde{\vk{w}}_v(\tau_1,\tau_2) := \widetilde{\Sigma}_v^{-1}(\tau_1,\tau_2)\vk{b}_v(\tau_1,\tau_2).
\EQNY
By conditioning on $\widetilde{\vk{X}}_v(\tau_1,\tau_2) = v\vk{b}_v(\tau_1,\tau_2)+\vk{x}/v$ and using the law of total probability, we obtain, continuing \eqref{eq:p_tau_v},
\BQN
&&p(\tau_1,\tau_2;v)\nonumber \\ 
&& \le v^{-\abs{I}} \int_{\R^{\abs{I}}}
\varphi_{\widetilde{\Sigma}_v(\tau_1,\tau_2)}\left(
v\vk{b}_v(\tau_1,\tau_2)-\frac{\vk{x}}{v}\right)
\pk{\exists t_1 \in [\tau_1,\tau_1+1],t_2 \in [\tau_2,\tau_2+1]:\vk{\chi}_v(t_1,t_2)>\vk{x}}
d\vk{x} \nonumber
\\
&& \leq v^{-\abs{I}} \varphi_{\widetilde{\Sigma}_v(\tau_1,\tau_2)}\left(
v\vk{b}_v(\tau_1,\tau_2)\right)
\int_{\R^{\abs{I}}} e^{(\widetilde{\vk{w}}_v(\tau_1,\tau_2))^{\top} \vk{x}}
\pk{\exists t_1 \in [\tau_1,\tau_1+1],t_2 \in [\tau_2,\tau_2+1]:\vk{\chi}_{v}(t_1,t_2)>\vk{x}}
d\vk{x}, \label{eq:p_tau_v-2}
\EQN
where
 $$\vk{\chi}_v(t_1,t_2) := v\left(\widetilde{\vk{X}}_v(t_1,t_2) - v\vk{b}_v(t_1,t_2)\right)+\vk{x}\Big{|} \left(
v\left(\widetilde{\vk{X}}_v(\tau_1,\tau_2) - v\vk{b}_v(\tau_1,\tau_2)\right)+\vk{x}=0
\right),$$
and in the last inequality we used the inequality
$$
\varphi_{\widetilde{\Sigma}_v(\tau_1,\tau_2)}\left(
v\vk{b}_v(\tau_1,\tau_2)-\frac{\vk{x}}{v}\right) \le   \varphi_{\widetilde{\Sigma}_v(\tau_1,\tau_2)}\left(
v\vk{b}_v(\tau_1,\tau_2)\right) \cdot   e^{(\widetilde{\vk{w}}_v(\tau_1,\tau_2))^{\top} \vk{x}}.
$$

\COM{
$\vk{X}_v(t)\coloneq \left(
\frac{A \vk{B}_H(t_0+tv^{-1/H})}{(t_0+tv^{-1/H})^H}
\right)_I$, $\widetilde{\vk{X}}_v(t,s)\coloneq
\frac{\vk{X}_v(t)+\vk{X}_v(s)}{2}$, $R_v(t,s) \coloneq \E{\vk{X}_v(t)\vk{X}^{\top}_u(s)}$, $\widetilde{R}_v(t_1,t_2,s_1,s_2) \coloneq \E{\widetilde{\vk{X}}_v(t_1,t_2)\widetilde{\vk{X}}^{\top}_v(s_1,s_2)}$, $\Sigma_v(t) \coloneq R_v(t,t)$, $\widetilde{\Sigma}_v(t_1,t_2) \coloneq \widetilde{R}(t_1,t_2,t_1,t_2)$, $\Sigma \coloneq A A^{\top}$, $$r_v(t,s) \coloneq
\E{
	\frac{B_{H,1}(t_0+tv^{-1/H})}{|t_0+tv^{-1/H}|^H}
	\frac{B_{H,1}(t_0+sv^{-1/H})}{|t_0+sv^{-1/H}|^H}
}
=
\frac{|t_0+tv^{-1/H}|^{2H}+|t_0+sv^{-1/H}|^{2H}-|t-s|^{2H} v^{-2}}{2|t_0+tv^{-1/H}|^{H}|t_0+sv^{-1/H}|^{H}},$$ $\sigma_v^2(t) \coloneq r_v(t,t)$, $\widetilde{r}_v(t_1,t_2,s_1,s_2) \coloneq \frac{1}{4} \left(r_v(t_1,s_1)+r_v(t_1,s_2)+r_v(t_2,s_1)+r_v(t_2,s_2)\right)
$.

} 
In the following, we shall derive suitable bounds for the \ccP{integral} and the term $ \varphi_{\widetilde{\Sigma}_v(\tau_1,\tau_2)}\left(v\vk{b}_v(\tau_1,\tau_2)\right)$ in \eqref{eq:p_tau_v-2}, respectively.
We start with the integral term, for which we shall apply \cite[Lemma 8]{Iev24}. 

 First, note that for any small  $\epsilon>0$ there exists $C_{M,\epsilon}>0$ and $v_{M,\epsilon}>0$ such that for all $v \geq v_{M,\epsilon}$ and all $t_1,s_1 \in [\tau_1,\tau_1 +1]$ and $t_2,s_2 \in [\tau_2,\tau_2+1]$ it holds that
\BQN
\label{rtilde_uasymphsmall}
v^2 \left|\widetilde{r}_v(t_1,t_2,s_1,s_2)-\left(1-\frac{\rJ{(\tau_2-\tau_1)}^{2H}v^{-2}}{4t_0^{2H}}\right)\right| &\leq& C_{M,\epsilon} +\epsilon (\tau_2-\tau_1)^{2H}.
\EQN
Indeed, by the Taylor's formula and the inequalities $|a^{2H}-b^{2H}| \leq |a-b|^{2H}$ for $a,b>0,H \in (0,1/2)$ and $|a^{2H}-b^{2H}| \leq 2H\max(a,b)^{2H-1} |a-b| \leq 2H( \epsilon \max(a,b)^{2H}+\epsilon^{-2H+1}|a-b|^{2H})$ for $a,b>0, H \in (1/2,1)$, we derive that
\BQN
\label{rtilde_uasymphsmall_prelim}
v^2 \left|r_v(t_i,s_j)-\left(1-\frac{|\tau_i-\tau_j|^{2H}v^{-2}}{2t_0^{2H}}\right)\right|  &\leq& C_{M,\epsilon} (|t_i-s_j|^2 v^{2-2/H})+
\left|
\frac{|t_i-s_j|^{2H} }{2|t_0+t_i v^{-1/H}|^{H}|t_0+s_i v^{-1/H}|^{H}}
-
\frac{|\tau_i-\tau_j|^{2H}}{2 t_0^{2H}}
\right|
\notag
\\
&\leq&
C_{M,\epsilon}(1+\left||t_i-s_j|^{2H}-|\tau_i-\tau_j|^{2H}\right|)+\epsilon |\tau_i-\tau_j|^{2H}
\notag
\\
&\leq&
C_{M,\epsilon}+\epsilon |\tau_i-\tau_j|^{2H}
\EQN
holds for any $ i,j \in\{1, 2\}$. Thus, by summing the inequalities \eqref{rtilde_uasymphsmall_prelim} over $1 \leq i,j \leq 2$, we  establish  \eqref{rtilde_uasymphsmall}.

Next, from \eqref{rtilde_uasymphsmall} it directly follows that there exists $C_{M,\epsilon}>0$ and $v_{M,\epsilon}>0$ such that for all $v \geq v_{M,\epsilon}$ and all $t_1,s_1, \widetilde{t}_1,\widetilde{s}_1 \in [\tau_1,\tau_1 +1]$ and $t_2,s_2, \widetilde{t}_2,\widetilde{s}_2 \in [\tau_2,\tau_2+1]$, 
\BQN
\label{rtilde_uasymphsmall2}
v^2 \left|\widetilde{r}_v(t_1,t_2,s_1,s_2)-\rJ{\widetilde{r}_v(\widetilde{t}_1,\widetilde{t}_2, \widetilde{s}_1,\widetilde{s}_2)}\right| \leq C_{M,\epsilon} +\epsilon (\tau_2-\tau_1)^{2H}.
\EQN
Additionally, there exists $C_M>0$ and $v_M>0$ such that for all $v \geq v_M$ and all $t,s \in [\tau_k,\tau_k+1]$ with $k \in \{1,2\}$, 
\BQN
\label{r_vhoelderhsmall}
v^2(1-r_v(t,s))\leq C_M|t-s|^{2H}.
\EQN
\COM{
Next, $R_v(t,s)=r_v(t,s) \Sigma$ and $\widetilde{R}_v(t_1,t_2,s_1,s_2) = \widetilde{r}_v(t_1,t_2,s_1,s_2)\Sigma$.

By the law of total probability, setting $\vk{\zeta}(t)\coloneq \frac{\vk{\mu}+\vk{\nu}t}{t^H}$, $\vk{b}_v(t_1,t_2) \coloneq
\frac{1}{2}
\left(
\vk{\zeta}_I(t_0+\frac{t_1}{2}v^{-1/H})
+
\vk{\zeta}_I(t_0+\frac{t_2}{2}v^{-1/H})
\right)
$, $\widetilde{\vk{w}} \coloneq \widetilde{\Sigma}_v^{-1}(\tau_1,\tau_2)\vk{b}_v(t_1,t_2)$ ($\widetilde{\vk{w}}$ may depend on $u,\tau_1,\tau_2$) and $$\vk{\chi}_v(t_1,t_2) \coloneq \left(v\left(\widetilde{\vk{X}}_v(t_1,t_2) - v\vk{b}_v(t_1,t_2)\right)+\vk{x} \mid
v\left(\widetilde{\vk{X}}_v(\tau_1,\tau_2) - v\vk{b}_v(\tau_1,\tau_2)\right)+\vk{x}=0
\right),$$ one obtains
\BQNY
&&p(\tau_1,\tau_2,1;v) \\&& \leq \pk{\exists t_1 \in [\tau_1,\tau_1+1],t_2 \in [\tau_2,\tau_2+1]:\widetilde{\vk{X}}_v(t_1,t_2)>v\vk{b}_v(\tau_1,\tau_2)}\\
&& =v^{-|I|} \int_{\R^I}
\varphi_{\widetilde{\Sigma}_v(\tau_1,\tau_2)}\left(
v\vk{b}_v(\tau_1,\tau_2)-\frac{\vk{x}}{v}\right)
\pk{\exists t_1 \in [\tau_1,\tau_1+1],t_2 \in [\tau_2,\tau_2+1]:\vk{\chi}_v(t_1,t_2)>\vk{x}}
d\vk{x}
\\
&& \leq v^{-|I|} \varphi_{\widetilde{\Sigma}_v(\tau_1,\tau_2)}\left(
v\vk{b}_v(\tau_1,\tau_2)\right)
\int_{\R^I} e^{\widetilde{\vk{w}}^{\top} \vk{x}}
\pk{\exists t_1 \in [\tau_1,\tau_1+1],t_2 \in [\tau_2,\tau_2+1]:\vk{\chi}_{v}(t_1,t_2)>\vk{x}}
d\vk{x}.
\EQNY
} 
 Note that, as $v\to \IF$,
$
\vk {\widetilde w}_v(\tau_1, \tau_2) \to \vk w/ t_0^{H} >\vk 0.
$
Thus, 
\BQN\label{eq:www}
\vk {\widetilde w}_v(\tau_1, \tau_2) \le  \frac{2}{t_0^{H}} \vk w =: \overline{ \vk w}
\EQN
holds for all large $v$.

\kk{In order to get a proper upper bound for \ccP{the integral term in} (\ref{eq:p_tau_v-2}), we shall apply \cite[Lemma 8]{Iev24},
for which we check the following three conditions
(recall that $F\subseteq I=\{1,\ldots,d\}$ is defined such that $\vk x_F>\vk 0$ and $\vk x_{I\setminus F}<\vk 0$ in the integral in \eqref{eq:p_tau_v-2}):}
\BQN
&&\sup_{F\subseteq I} \sup_{t,s\in[0,1]} \overline{\vk w}_F ^\top \E{\vk \chi_{v,F}(\tau_1+t,\tau_2+s)} \le C_{M,\epsilon}+ \epsilon (\tau_2-\tau_1)^{2H}+\epsilon \sum_{j=1}^d \abs{x_j},\label{eq:Pavel-1}\\
&&\sup_{F\subseteq I} \sup_{t,s\in[0,1]} \text{Var} \LT(\overline{\vk w}_F ^\top \vk \chi_{v,F}(\tau_1+t,\tau_2+s)\RT) \le C_{M,\epsilon}+\epsilon  (\tau_2-\tau_1)^{2H}, \label{eq:Pavel-2} 
\EQN
and, for any $F\subseteq I$ and $t_1,s_1 \in [\tau_1,\tau_1+1]$, $t_2,s_2 \in [\tau_2,\tau_2 + 1]$
\BQN
 \text{Var} \LT( \overline{\vk w}_F^\top \vk \chi_{v,F}(t_1,t_2) -\overline{\vk w}_F^\top \vk \chi_{v,F}(s_1,s_2)\RT)\le  C_M\LT( (|t_1-s_1|^{2H}+|t_2-s_2|^{2H} \RT) \label{eq:Pavel-3}
\EQN
hold for all  large enough $v$.

\kk{\underline{Inequality \eqref{eq:Pavel-1}}}.
Note that by Lagrange's mean value theorem there exist $C_M>0$, $v_M>0$ such that for $v \geq v_M$ and $t_1 \in [\tau_1,\tau_1+1]$ and $t_2 \in [\tau_2,\tau_2+1]$ it holds that, if $H < {1}/{2}$, then
$$v^2\ABs{\vk{b}_v(t_1,t_2)-\vk{b}_v(\tau_1,\tau_2)}
\leq v^{2-1/H} \sup_{s \in [\tau_1,\tau_1+1] \cup [\tau_2,\tau_2+1]} \ABs{\vk{\zeta}'_I(t_0+sv^{-1/H})}
\leq C_M,$$
and, if $H>1/2$ and $H \vk \nu_I = (1-H) t_0 \vk \mu_I$ (i.e., $\vk{\zeta}'_I(t_0)=0$), then, by Taylor's formula with Lagrange's remainder, 
\BQNY
v^2\Abs{\vk{b}_v(t_1,t_2)-\vk{b}_v(\tau_1,\tau_2)}
&\leq& v^2\Abs{\vk{b}_v(t_1,t_2)-\vk{\zeta}_I(t_0)}+v^2\Abs{\vk{b}_v(\tau_1,\tau_2)-\vk{\zeta}_I(t_0)}\\
&\leq&
v^{2-2/H} (|\tau_1|+|\tau_2|+2)^2
\sup_{|s| \leq \max(|\tau_1|,|\tau_2|)+1} |\vk{\zeta}''_I(t_0+sv^{-1/H})|
\\
&\leq&
C_M.
\EQNY
Hence,  for any small $\epsilon>0$ there exist $C_{M,\epsilon},v_{M,\epsilon}>0$ such that for all $v \geq v_{M,\epsilon}$ and all $t_1 \in [\tau_1,\tau_1+1], t_2 \in [\tau_2,\tau_2+1]$, 
\BQNY
&&\left|\E{\vk{\chi}_v(t_1,t_2)}\right|
\notag
=
\left|\left(-v^2 \vk{b}_v(t_1,t_2) + \vk{x}\right) - \widetilde{R}_v(t_1,t_2,\tau_1,\tau_2)\widetilde{\Sigma}_v^{-1}(\tau_1,\tau_2)
\left(-v^2 \vk{b}_v(\tau_1,\tau_2)+\vk{x}\right)\right|
\notag
\\&&\ \ \leq
v^2
\left|
\left(\widetilde{\Sigma}_v(\tau_1,\tau_2)-\widetilde{R}_v(t_1,t_2,\tau_1,\tau_2)\right)
\widetilde{\Sigma}_v^{-1}(\tau_1,\tau_2)
\vk{b}_v(\tau_1,\tau_2)
\right|+C_{M,\epsilon}+\epsilon \left|\vk{x}\right|
\notag
\\&&
\ \ \leq C_{M,\epsilon} v^2\left|\left(\widetilde{r}_v(\tau_1,\tau_2,\tau_1,\tau_2)-\widetilde{r}_v(t_1,t_2,\tau_1,\tau_2)\right)
\right| + C_{M,\epsilon}+\epsilon \left|\vk{x}\right|
\notag
\\&&
\ \ \leq
C_{M,\epsilon}+ \epsilon ((\tau_2-\tau_1)^{2H}+\left|\vk{x}\right|),
\EQNY
where the second inequality follows by the fact that $\lim_{v\to\IF} \widetilde{r}_v^{-1}(\tau_1,\tau_2,\tau_1,\tau_2)  \vk{b}_v(\tau_1,\tau_2) = \vk b/ t_0^{H}$ and the third (last) inequality follows  by using \eqref{rtilde_uasymphsmall2}. This yields \eqref{eq:Pavel-1}.

\kk{\underline{Inequality \eqref{eq:Pavel-2}}}. Note that, for any $t_1 \in [\tau_1,\tau_1+1], t_2 \in [\tau_2,\tau_2+1]$, 
\BQNY
\mathrm{Var}(\overline{\vk{w}}_F^{\top} \vk{\chi}_{v,F}(t_1,t_2)) = \overline{\vk{w}}_F^{\top} \Sigma_{FF} \overline{\vk{w}}_F
\cdot K_v(t_1,t_2,t_1,t_2),
\EQNY
where  $K_v(t_1,t_2,s_1,s_2) :=
v^2
\left( \widetilde{r}_v(t_1,t_2,s_1,s_2)-\widetilde{r}_v(t_1,t_2,\tau_1,\tau_2)\widetilde{r}_v^{-1}(\tau_1,\tau_2,\tau_1,\tau_2)\widetilde{r}_v(\tau_1,\tau_2,s_1,s_2)
\right)
$. It follows from \eqref{rtilde_uasymphsmall2} that, for any small $\epsilon>0$ there exist $C_{M,\epsilon},v_{M,\epsilon}>0$ such that for all $v \geq v_{M,\epsilon}$ and all $t_1,s_1 \in [\tau_1,\tau_1+1], t_2,s_2 \in [\tau_2,\tau_2+1]$,
\BQNY
|K_v(t_1,t_2,s_1,s_2)|
&\leq& C_{M,\epsilon}v^2\ABs{\widetilde{r}_v(t_1,t_2,s_1,s_2)\widetilde{r}_v(\tau_1,\tau_2,\tau_1,\tau_2)-\widetilde{r}_v(t_1,t_2,\tau_1,\tau_2)\widetilde{r}_v(\tau_1,\tau_2,s_1,s_2)}
\notag
\\&\leq&
C_{M,\epsilon}+\epsilon |\tau_1-\tau_2|^{2H}.
\EQNY
This implies \eqref{eq:Pavel-2}.

\kk{\underline{Inequality  \eqref{eq:Pavel-3}}}.
We have, using \eqref{cond}, that for $t_1,s_1 \in [\tau_1,\tau_1+1]$, $t_2,s_2 \in [\tau_2,\tau_2 + 1]$,
\BQNY
\Var\left(\overline{\vk{w}}_F^{\top} \vk{\chi}_{v,F}(t_1,t_2)-\overline{\vk{w}}_F^{\top} \vk{\chi}_{v,F}(s_1,s_2)\right)
\leq v^2 \Var \left(\overline{\vk{w}}_F^{\top} \widetilde{\vk{X}}_{v,F}(t_1,t_2)-\overline{\vk{w}}_F^{\top} \widetilde{\vk{X}}_{v,F}(s_1,s_2)\right).
\EQNY
Furthermore, there exist $C_M>0$, $v_M>0$ such that for all $v \geq v_M$ and all
$t_1,s_1 \in [\tau_1,\tau_1+1]$, $t_2,s_2 \in [\tau_2,\tau_2 + 1]$, 
\BQN
\label{Pavel_Hoelder}
&&
v^2 \Var \left(\overline{\vk{w}}_F^{\top} \widetilde{\vk{X}}_{v,F}(t_1,t_2)-\overline{\vk{w}}_F^{\top} \widetilde{\vk{X}}_{v,F}(s_1,s_2)\right)
\notag
\\
&&
\leq
\frac{v^2}{2}\mathrm{Var}\left(
\overline{\vk w}^{\top}_F \widetilde{\vk{X}}_{v,F}(t_1)-\overline{\vk w}^{\top}_F \widetilde{\vk{X}}_{v,F}(s_1)
\right)
+
\frac{v^2}{2}\mathrm{Var}\left(
\overline{\vk w}^{\top}_F \widetilde{\vk{X}}_{v,F}(t_2)-\overline{\vk w}^{\top}_F \widetilde{\vk{X}}_{v,F}(s_2)
\right)
\notag
\\
&&
\leq C_M (|t_1-s_1|^{2H}+|t_2-s_2|^{2H}),
\EQN
where the last inequality follows from \eqref{r_vhoelderhsmall}. Thus,  \eqref{eq:Pavel-3} is established.
\COM{
Moreover, since the conditional variance is always less than or equal to the unconditional one, we have, for any $F\subseteq \{1,\ldots,d\}$, that
\BQNY
&& \text{Var} \LT( \vk w_F (\vn_2)^\top \vk \chi_{v,F}(\tau_1+t,\tau_2+s) -\vk w_F (\vn_2)^\top \vk \chi_{v,F}(\tau_1+t_1,\tau_2+s_1)\RT)\\
&\le &  v^2   \text{Var} \LT( \vk w_F (\vn_2)^\top      \widetilde{\vk{X}}_{v,F}(\tau_1+t,\tau_2+s) -\vk w_F (\vn_2)^\top     \widetilde{\vk{X}}_{v,F}(\tau_1+t_1,\tau_2+s_1)\RT)\\
&\le&\frac{v^2}{2}  \text{Var} \LT( \vk w_F (\vn_2)^\top      \vk{X}_{v,F}(\tau_1+t) -\vk w_F (\vn_2)^\top     \vk{X}_{v,F}(\tau_1+t_1)\RT)\\
&&+ \frac{v^2}{2}  \text{Var} \LT( \vk w_F (\vn_2)^\top      \vk{X}_{v,F}(\tau_1+s) -\vk w_F (\vn_2)^\top     \vk{X}_{v,F}(\tau_1+s_1)\RT)\\
&\le & G_2\LT( (|t-s|^{2H}+|t_1-s_1|^{2H} \RT),
\EQNY
where $G_2>0$ is a constant which can be chosen independent of $F$. Consequently, by applying Lemma ?? in [Pavel's paper], we obtain, for the integral in \eqref{eq:pttTv}, that 
\BQN\label{eq:BoundInt}
\int_{\R^d} e^{\LT(\vk {\widetilde w}_v(\tau_1, \tau_2) \RT)^\top \vk x}\,  \pk{\exists t \in [\tau_1, \tau_1+ T], s \in [\tau_2,\tau_2+T]: \vk \chi_v(t,s)  > \vk x} d\vk x \le C_0 e^{G_0 \vn_1 (\tau_2-\tau_1)^{2H}}
\EQN
holds for all large enough $v$, where $C_0, G_0>0$ are some constants.

 Remark that by Lagrange's mean value theorem there exist $C_M>0$, $v_M>0$ such that for $v \geq v_M$ and $t_1 \in [\tau_1,\tau_1+1]$ and $t_2 \in [\tau_2,\tau_2+1]$ it holds that
$$v^2\ABs{\vk{b}_v(t_1,t_2)-\vk{b}_v(\tau_1,\tau_2)}
\leq v^{2-1/H} \sup_{s \in [\tau_1,\tau_1+1] \cup [\tau_2,\tau_2+1]} \ABs{\vk{\zeta}'_I(t_0+sv^{-1/H})}
\leq C_M,$$
if $H < \frac{1}{2}$. On the other hand, if $H>1/2$, since in this case $\vk{\zeta}'_I(t_0)=0$, by Taylor's formula with Lagrange's remainder one obtains
\BQNY
v^2\Abs{\vk{b}_v(t_1,t_2)-\vk{b}_v(\tau_1,\tau_2)}
&\leq& v^2\Abs{\vk{b}_v(t_1,t_2)-\zeta_I(t_0)}+v^2\Abs{\vk{b}_v(\tau_1,\tau_2)-\zeta_I(t_0)}\\
&\leq&
v^{2-2/H} (|\tau_1|+|\tau_2|+2)^2
\sup_{|s| \leq \max(|\tau_1|,|\tau_2|)+1} |\vk{\zeta}''_I(t_0+sv^{-1/H})|
\\
&\leq&
C_M.
\EQNY
Hence, by \eqref{rtilde_uasymphsmall2} for each $\epsilon>0$ there exist $C_{M,\epsilon},v_{M,\epsilon}>0$ such that if $v \geq v_{M,\epsilon}$ and $t_1 \in [\tau_1,\tau_1+1], t_2 \in [\tau_2,\tau_2+1]$, then
\BQN
\label{Pavel_mean}
&&\left|\E{\vk{\chi}_v(t_1,t_2)}\right|
\notag
=
\left|\left(-v^2 \vk{b}_v(t_1,t_2) + \vk{x}\right) - \widetilde{R}_v(t_1,t_2,\tau_1,\tau_2)\widetilde{\Sigma}_v^{-1}(\tau_1,\tau_2)
\left(-v^2 \vk{b}_v(\tau_1,\tau_2)+\vk{x}\right)\right|
\notag
\\&&\leq
v^2
\left|
\left(\widetilde{\Sigma}_v(\tau_1,\tau_2)-\widetilde{R}_v(t_1,t_2,\tau_1,\tau_2)\right)
\widetilde{\Sigma}_v^{-1}(\tau_1,\tau_2)
\vk{b}_v(\tau_1,\tau_2)
\right|+C_{M,\epsilon}+\epsilon \left|\vk{x}\right|
\notag
\\&&
\leq C_{M,\epsilon}\left|\left(\widetilde{r}_v(\tau_1,\tau_2,\tau_1,\tau_2)-\widetilde{r}_v(t_1,t_2,\tau_1,\tau_2)\right)
\widetilde{r}_v^{-1}(\tau_1,\tau_2,\tau_1,\tau_2)
\left(-v^2 \vk{b} - \vk{\mu}_I v^{2-1/H} \frac{\tau_1+\tau_2}{2}\right)\right|
\notag
\\
&&\quad
+\epsilon(|\tau_2-\tau_1|^{2H}+\left|\vk{x}\right|)+C_{M,\epsilon}
\notag
\\&&
\leq
C_{M,\epsilon}+ \epsilon(|\tau_2-\tau_1|^{2H})+\left|\vk{x}\right|).
\EQN
By \eqref{cond}
\BQNY
\Var\left(\widetilde{\vk{w}}_F^{\top} \vk{\chi}_{v,F}(t_1,t_2)-\widetilde{\vk{w}}_F^{\top} \vk{\chi}_{v,F}(s_1,s_2)\right)
\leq v^2 \Var \left(\widetilde{\vk{w}}_F^{\top} \widetilde{\vk{X}}_{v,F}(t_1,t_2)-\widetilde{\vk{w}}_F^{\top} \widetilde{\vk{X}}_{v,F}(s_1,s_2)\right).
\EQNY
Next, as for $v$ large enough $\widetilde{\vk{w}}$ is uniformly bounded in $v, \tau_1,\tau_2$, there exist $C_M>0$, $v_M>0$ such that for each $v \geq v_M$,
$t_1,s_1 \in [\tau_1,\tau_1+1]$, $t_2,s_2 \in [\tau_2,\tau_2 + 1]$ it holds that
\BQN
\label{Pavel_Hoelder}
&&
v^2 \Var \left(\widetilde{\vk{w}}_F^{\top} \widetilde{\vk{X}}_{v,F}(t_1,t_2)-\widetilde{\vk{w}}_F^{\top} \widetilde{\vk{X}}_{v,F}(s_1,s_2)\right)
\notag
\\
&&
\leq
\frac{v^2}{2}\mathrm{Var}\left(
\widetilde{w}^{\top}_F \widetilde{\vk{X}}_{v,F}(t_1)-\widetilde{w}^{\top}_F \widetilde{\vk{X}}_{v,F}(s_1)
\right)
+
\frac{v^2}{2}\mathrm{Var}\left(
\widetilde{w}^{\top}_F \widetilde{\vk{X}}_{v,F}(t_2)-\widetilde{w}^{\top}_F \widetilde{\vk{X}}_{v,F}(s_2)
\right)
\notag
\\
&&
\leq C_M (|t_1-s_1|^{2H}+|t_2-s_2|^{2H}),
\EQN
where the last inequality follows from \eqref{r_vhoelderhsmall}.
} 

Consequently, an application of \cite[Lemma 8]{Iev24} yields that, for any small $\epsilon>0$ there exist $C_{M,\epsilon},v_{M,\epsilon}>0$ such that, for all $v \geq v_{M,\epsilon}$,
\BQN
\int_{\R^{\abs{I}}} e^{(\widetilde{\vk{w}}_v(\tau_1,\tau_2))^{\top} \vk{x}}
\pk{\exists t_1 \in [\tau_1,\tau_1+1],t_2 \in [\tau_2,\tau_2+1]:\vk{\chi}_{v}(t_1,t_2)>\vk{x}}
d\vk{x} \leq e^{C_{M,\epsilon}+\epsilon(\tau_2-\tau_1)^{2H}}.\label{eq:Int_w}
\EQN
It remains to estimate $\varphi_{\widetilde{\Sigma}_v(\tau_1,\tau_2)}\left(
v\vk{b}_v(\tau_1,\tau_2)\right)$. By  Taylor's formula with Lagrange's remainder we derive that there exist $C_{M},v_{M}>0$ such that, for all $v \geq v_{M}$,
\BQNY
v^2\ABs{\vk{b}_v(\tau_1,\tau_2)-\vk{\zeta}_I\left(t_0+\frac{\tau_1+\tau_2}{2}v^{-1/H}\right)}
\leq v^{2-2/H}|\tau_1-\tau_2|^{2}
\sup_{s \in [\tau_1,\tau_2]}
|\vk{\zeta}''_I(t_0+sv^{-1/H})|
\leq C_M.
\EQNY
Hence, there exists $C>0$ such that for any small $\epsilon>0$ there exist $C_{M,\epsilon}, C_{M}, v_{M,\epsilon}>0$ such that, for all $v \geq v_{M,\epsilon}$,
\BQNY
\varphi_{\widetilde{\Sigma}_v(\tau_1,\tau_2)}\left(
v\vk{b}_v(\tau_1,\tau_2)\right)
&\leq&C_{M,\epsilon}\exp\left(-\frac{v^2}{2}
\vk{b}_v^{\top}(\tau_1,\tau_2)
\widetilde{\Sigma}_v^{-1}(\tau_1,\tau_2)
\vk{b}_v(\tau_1,\tau_2)
\right)
\\&=&
C_{M,\epsilon}\exp\left(-\frac{v^2}{2}\widetilde{r}_v^{-1}(\tau_1,\tau_2,\tau_1,\tau_2)
\vk{b}_v^{\top}(\tau_1,\tau_2)
\Sigma_{II}^{-1}\vk{b}_v(\tau_1,\tau_2)
\right)
\\&\leq&
C_{M,\epsilon}\exp\left(-\frac{v^2}{2}\widetilde{r}_v^{-1}(\tau_1,\tau_2,\tau_1,\tau_2) g_I\left(t_0+\frac{\tau_1+\tau_2}{2}v^{-1/H}
\right)
\right)
\\&\leq&
C_{M,\epsilon}
\exp\left(-\left(\frac{v^2}{2}
+\left(\frac{1}{4t_0^{2H}}-\epsilon\right)(\tau_2-\tau_1)^{2H}
\right)
g_I\left(t_0+\frac{\tau_1+\tau_2}{2}v^{-1/H}
\right)
\right)
\\&\leq&
C_{M,\epsilon}
\exp\left(-\frac{v^2}{2} g_I\left(t_0\right)
-C_{M}(\tau_2-\tau_1)^{2H}
\right),
\EQNY
where the penultimate inequality follows from \eqref{rtilde_uasymphsmall}. This, together with \eqref{eq:p_tau_v-2} and \eqref{eq:Int_w}, establishes the claim of Lemma \ref{Hsmalldoubleevent}. 
\QED


{\bf Proof of Lemma \ref{lem:Pick2}:}
The proof is analogous to the proof of Lemma \ref{lem:Pick1}, and thus we shall only present the main differences and some key calculations. We define
\BQNY
&&
\vk{X}_{v,\tau}(t) := \vk{X}(t_0+\tau v^{-2} + t v^{-2}), \ \ \Sigma_{v,\tau} = \E{\vk{X}_{v,\tau}(0) \vk{X}_{v,\tau}(0)^{\top}}=(t_0+\tau v^{-2})^{2H}\Sigma,
\\
&&
\vk{Z}_{v,\tau}(t)
:=
\overline{\vk{v}}\left(
\vk{X}_{v,\tau}(t)-\widetilde{\vk{b}}v-\tau \vk{\mu} v^{-1}-t\vk{\mu}v^{-1}
\right)
+\vk{x},
\\
&&
r_{v,\tau}(t,s):=\frac{1}{2}\LT( (t_0+\tau v^{-2}+tv^{-2})^{2H} + (t_0+\tau v^{-2}+sv^{-2})^{2H} -\abs{t-s}^{2H} v^{-4H}\RT).
\EQNY
Then
\BQNY
&&
\pk{\exists_{t \in[t_0+\tau v^{-2}, t_0+(\tau+T) v^{-2}]}
	\vk X(t) > (\vk \nu +\vk \mu t) v }
\\
&&
= v^{-\abs{I}} \int_{\R^d}  \pk{\exists_{t \in[0, T]}
	\vk Z_{v,\tau}(t) > (\vk b- \wtb) \olv v+\vk x \mid  \vk Z_{v,\tau}(0)=\vk 0} \varphi_{\Sigma_{v,\tau}}(\wtb  v+ \tau \vk \mu v^{-1}-\vk x/\olv) d\vk x.
\EQNY
Define
\BQNY
\vk \chi_{v,\tau}(t) := (\vk Z_{v,\tau}(t) \mid \vk Z_{v,\tau}(0)=\vk 0).
\EQNY
Similarly as in the proof of Lemma \ref{lem:Pick1}, we derive that
\BQNY
\E{\vk \chi_{v,\tau}(t) }
= \LT(1-  r_{v,\tau}(t,0) r_{v,\tau}(0,0)^{-1}\RT) \vk x + \overline{\vk v}\LT( ( r_{v,\tau}(t,0) r_{v,\tau}(0,0)^{-1} -1) ( \wtb v+ \tau \vk \mu v^{-1} )  -t \vk \mu v^{-1}\RT).
\EQNY
Next, it follows that
\BQNY
r_{v,\tau}(t,0) r_{v,\tau}(0,0)^{-1} -1 = \frac{(t_0+\tau v^{-2}+tv^{-2})^{2H} - (t_0+\tau v^{-2})^{2H}}{2 (t_0+\tau v^{-2})^{2H}} -  \frac{t^{2H} v^{-4H}}{2 (t_0+\tau v^{-2})^{2H}},
\EQNY
and
\BQNY
(t_0+\tau v^{-2}+tv^{-2})^{2H} - (t_0+\tau v^{-2})^{2H} = 2H (t_0+\tau v^{-2})^{2H-1} tv^{-2} + H(2H-1) (t_0+\tau' v^{-2})^{2H-2} (tv^{-2})^2, 
\EQNY
with some $\tau' \in [\tau, \tau+t]$. Some elementary calculations yield that, as $v\to\IF,$
\BQNY
\E{\vk \chi_{v,\tau}(t) }
\to  \begin{pmatrix}
	\LT(\frac{H}{t_0}\vk b_I - \vk{\mu}_I\RT) t \\
	\vk 0 _{I^c}
\end{pmatrix}
\EQNY
holds uniformly for  $\tau$ such that $\abs{\tau} \leq T(\hN_v+1)$.
Similarly to the calculations for the mean, we obtain that, as $v\to\IF,$
\BQNY
\E{[\vk \chi_{v,\tau}(t) -\E{\vk \chi_{v,\tau}(t)}][\vk \chi_{v,\tau}(s) -\E{\vk \chi_{v,\tau}(s)}]^\top}\to  \vk{0}
\EQNY
holds uniformly for  $\tau$ such that $\abs{\tau} \leq T(\hN_v+1)$.
Thus, 
as $v\to\IF,$
\BQNY
\pk{\exists_{t \in[0, T]}
	\vk Z_{v,\tau}(t) > (\vk b- \wtb) \olv v +\vk x \mid  \vk Z_{v,\tau}(0)=\vk 0}  \to    \pk{\exists_{t \in[0, T]}
	\left(\frac{H}{ t_0}\vk{b}_I-\vk{\mu}_I\right)t
	>\vk x_I }\cdot   \kls{\mathbb {I}}_{(\vk x_K < \vk 0_K)}.
\EQNY
Proceeding similarly to the proof of
Lemma \ref{lem:Pick1}, we can obtain
\BQNY
&&\pk{\exists_{t \in[t_0+\tau v^{-2}, t_0+(\tau+T) v^{-2}]}
 \vk X(t) > (\vk \nu +\vk \mu t) v }\\
&&\sim  v^{-\abs{I}} \frac{\widetilde{\mathcal H}_{I}(T)}{\sqrt{ (2\pi t_0^{2H})^d \abs{\Sigma}}} e^{-\frac{v^2}{2} g_I(t_0+\tau v^{-2})}\\
&&\times \int_{\R^{\abs{I^c}}} e^{-\frac{1}{2 t_0^{2H}} \vk x_{I^c}^ \top (\Sigma^{-1})_{I^cI^c} \vk x_{I^c}}   \kls{\mathbb {I}}_{\LT(\vk x_K <  - ( \tau v^{-1}) [[(\Sigma^{-1})_{I^cI^c}]^{-1} (\Sigma^{-1} \vk \mu)_{I^c}]_K \RT)} d\vk x_{I^c},
\EQNY
where \BQNY
\widetilde{\mathcal H}_{I}(T)&=& \int_{\R^{\abs{I}}} e^{\frac{1}{t_0^{2H}} \vk w_I^\top \vk x_I}   \pk{\exists_{t \in[0, T]}
	\vk{x}_I<-\left(\vk{\mu}_I-\frac{H}{ t_0}\vk{b}_I\right)t
	} d\vk x_I,
\\
&=&
\frac{t_0^{2H|I|}}{\prod_{i \in I}w_i}
\int_{\R^{\abs{I}}} e^{ \vk 1_I^\top \vk x_I}   \pk{\exists_{t \in[0, T]}
	\vk{x}_I<-\frac{1}{t_0^{2H}} \mathrm{diag}(\vk{w}_I)\left(\vk{\mu}_I-\frac{H}{ t_0}\vk{b}_I\right)t
} d\vk x_I.
\EQNY
Furthermore, since $\vk w_I>\vk 0_I$ and
\BQNY
\vk{1}_I^{\top}
\mathrm{diag}(\vk{w}_I)\left(\vk{\mu}_I-\frac{H}{ t_0}\vk{b}_I\right)
=\vk{b}_I^{\top}\Sigma_{II}^{-1}
\left(
\vk{b}_I'(t_0)-\frac{H}{t_0}\vk{b}_I
\right)
=\frac{t_0^{2H}}{2}g_I'(t_0)=0,
\EQNY
we obtain, by \cite[Lemma 5.3]{DHW},
\BQNY
\widetilde{\mathcal{H}}_I (T)
=
\frac{t_0^{2H|I|}}{\prod_{i \in I}w_i}
\LT( 1+\frac{T}{t_0^{2H}}\sum_{i \in I}
\left(w_i \mu_i - \frac{H}{t_0} w_i b_i\right)_{-} \RT),
\EQNY
with $\sum_{i \in I}
\left(w_i \mu_i - \frac{H}{t_0} w_i b_i\right)_{-}>0.$
This \ksn{combined with \eqref{integr_expr}}
completes the proof. \QED


{\bf Proof of Lemma \ref{Hbigdoubleevent}:}
	The proof proceeds
	analogously to the proof of Lemma \ref{Hsmalldoubleevent},
	however, some important changes should be applied. Hereafter the exact values of the constants $\epsilon, C_M,v_M,C_{M,\epsilon}$ and $v_{M,\epsilon}$ are not important and can be changed  from line to line. Moreover, all the inequalities hold  uniformly in $\tau_1,\tau_2 \in [-Mv-1,Mv]$ such that $\tau_2 -\tau_1 \geq 1$ and the constants do not depend on $\tau_1,\tau_2$.

Without loss of generality, we shall assume  $I = \{1,\ldots, d\}$.
Let $V \subset I$ be a non-empty set to be chosen later ($V$ will not depend on $\tau_1,\tau_2$). Denote
	\BQNY && \widetilde{\vk{X}}_v(t_1,t_2) := \begin{pmatrix}
		\vk{X}_{v,V}(t_1)\\
		\vk{X}_{v,I \backslash V}(t_2)
	\end{pmatrix}, \ \ \vk{X}_v(t):=  \frac{\vk{X}(t_0+tv^{-2})}{(t_0+tv^{-2})^{H}},\\ 
	&&\widetilde{R}_v(t_1,t_2,s_1,s_2) := \E{\widetilde{\vk{X}}_v(t_1,t_2)\widetilde{\vk{X}}^{\top}_v(s_1,s_2)},\ \  \widetilde{\Sigma}_v(t_1,t_2) := \widetilde{R}_v(t_1,t_2,t_1,t_2), 
	\\
	&&r_v(t,s) :=
	\E{
		\frac{B_{H,1}(t_0+tv^{-2})}{(t_0+tv^{-2})^H}
		\frac{B_{H,1}(t_0+sv^{-2})}{(t_0+sv^{-2})^H}
	}
	=
	\frac{(t_0+tv^{-2})^{2H}+(t_0+sv^{-2})^{2H}-|t-s|^{2H} v^{-4H}}{2 (t_0+tv^{-2})^{H}(t_0+sv^{-2})^{H}}.
	\EQNY
It follows that, for any $t \in [\tau_k,\tau_k +1]$ and $s \in [\tau_l,\tau_l+1]$, $k,l \in \{1,2\}$,
\BQNY
0 \leq
1-r_v(t,s)=\frac{\LT((t_0+tv^{-2})^{H}-(t_0+sv^{-2})^{H}\RT)^2 + (|t-s|v^{-2})^{2H}   }{2 (t_0+tv^{-2})^{H}(t_0+sv^{-2})^{H}}.
\EQNY
By Taylor's formula and the fact that $H>1/2$, we have that,
 for any small $\epsilon>0$
	there exists 
$v_{M,\epsilon}>0$ such that for all $v \geq v_{M,\epsilon}$ and all $t \in [\tau_k,\tau_k +1]$, $s \in [\tau_l,\tau_l+1]$, $k,l \in \{1,2\}$, 
\BQN
v^2\left(1-r_v(t,s)\right) &\leq&	\epsilon \abs{t-s}	\label{r_uhoelderhbig} \\
&\le &  \epsilon+\epsilon(\tau_2-\tau_1). \label{r_uasymphbig}
	\EQN

 Further, denote
\BQNY
\vk{b}_v(t_1,t_2):=
	\begin{pmatrix}
		\vk{\zeta}_V(t_0+t_1 v^{-2}) \\
		\vk{\zeta}_{I \backslash V}(t_0+t_2 v^{-2})
	\end{pmatrix}
\ \  \text{with} \ \ \vk{\zeta}(t) := \frac{\vk{\nu}+\vk{\mu}t}{t^H},
\EQNY
and
\BQNY
	&&\vk{\chi}_v(t_1,t_2) :=v\left(\widetilde{\vk{X}}_v(t_1,t_2) - v\vk{b}_v(t_1,t_2)\right)+\vk{x} \Big{|}  \left(
	v\left(\widetilde{\vk{X}}_v(\tau_1,\tau_2) -v \vk{b}_v(\tau_1,\tau_2)\right)+\vk{x}=0
	\right),\\
&& \widetilde{\vk{w}}_v(\tau_1, \tau_2) = \widetilde{\Sigma}_v^{-1}(\tau_1,\tau_2)\vk{b}_v(\tau_1,\tau_2).
\EQNY
		By the law of total probability, we obtain
	\BQN
	\hp(\tau_1,\tau_2;v) & \leq &\pk{\exists t_1 \in [\tau_1,\tau_1+1],t_2 \in [\tau_2,\tau_2+1]:\widetilde{\vk{X}}_v(t_1,t_2)>v \vk{b}_v(t_1,t_2)}\nonumber\\
	& \leq & v^{-|I|} \varphi_{\widetilde{\Sigma}_v(\tau_1,\tau_2)}\left(
	v\vk{b}_v(\tau_1,\tau_2)\right)
	\int_{\R^{\abs{I}}} e^{(\widetilde{\vk{w}}_v(\tau_1, \tau_2))^{\top} \vk{x}}
	\pk{\exists t_1 \in [\tau_1,\tau_1+1],t_2 \in [\tau_2,\tau_2+1]:\vk{\chi}_{v}(t_1,t_2)>\vk{x}}
	d\vk{x}.\nonumber\\ \label{int}
	\EQN
Note that, \kk{similarly to} \eqref{eq:www},
$\vk {\widetilde w}_v(\tau_1, \tau_2) \le   \overline{ \vk w}$
holds for all large $v$.

\kk{In order to find a tight  upper  bound for the integral in (\ref{int}), similarly to the proof of Lemma \ref{Hsmalldoubleevent}, we shall verify  conditions of  \cite[Lemma 8]{Iev24} as stated in
\eqref{eq:Pavel-1}-\eqref{eq:Pavel-3} \qq{with $H$ replaced by $1/2$}.}

\kk{\underline{Inequality \eqref{eq:Pavel-1}}}.
Recall that $F\subseteq I=\{1,\ldots,d\}$ is defined such that $\vk x_F>\vk 0$ and $\vk x_{I\setminus F}<\vk 0$ in the integral.
There exist $C_M>0$, $v_M>0$ such that for all $v \geq v_M$ and all $t_1 \in [\tau_1,\tau_1+1]$ and $t_2 \in [\tau_2,\tau_2+1]$, 
	$$v^2\ABs{\vk{b}_v(t_1,t_2)-\vk{b}_v(\tau_1,\tau_2)}
	\leq  \sup_{s \in [\tau_1,\tau_1+1] \cup [\tau_2,\tau_2+1]} \ABs{\vk{\zeta}'_I(t_0+sv^{-2})}
	\leq C_M,$$
	and therefore, by
	\eqref{r_uasymphbig}, for any small $\epsilon>0$
	there exist $C_{M,\epsilon}>0$ and $v_{M,\epsilon}>0$ such that for all $v \geq v_{M,\epsilon}$ and all $t_1 \in [\tau_1,\tau_1 +1]$, $t_2 \in [\tau_2,\tau_2+1]$,
\BQNY
&&\ABs{\E{\vk{\chi}_v(t_1,t_2)}}=
\ABs{\left(-v^2 \vk{b}_v(t_1,t_2) + \vk{x}\right) - \widetilde{R}_v(t_1,t_2,\tau_1,\tau_2)\widetilde{\Sigma}_v^{-1}(\tau_1,\tau_2)
	\left(-v^2 \vk{b}_v(\tau_1,\tau_2)+\vk{x}\right)}
\notag
\\&&\leq
v^2
\ABs{
	\left(\widetilde{\Sigma}_v(\tau_1,\tau_2)-\widetilde{R}_v(t_1,t_2,\tau_1,\tau_2)\right)
	\widetilde{\Sigma}_v^{-1}(\tau_1,\tau_2)
	\vk{b}_v(\tau_1,\tau_2)
}+C_{M,\epsilon}+\epsilon \left|\vk{x}\right|
\notag
\\&&
\le C_{M,\epsilon}+\epsilon((\tau_2-\tau_1)+\left|\vk{x}\right|),
\EQNY
where in the last inequality we use  (here $||A|| = \max_{i,j} \abs{a_{ij}}$ for a matrix $A$,  and \eqref{r_uasymphbig} is applied)
\BQNY
 \kk{v^2\left\Arrowvert\widetilde{\Sigma}_v(\tau_1,\tau_2)-\widetilde{R}_v(t_1,t_2,\tau_1,\tau_2) \right\Arrowvert}
 &=& \left\Arrowvert \begin{array}{cc}
		v^2(r_v(\tau_1,\tau_1)- r_v(t_1,\tau_1))\Sigma_{V,V} &
		v^2(r_v(\tau_1,\tau_2)- r_v(t_1,\tau_2)) \Sigma_{V,I \backslash V}\\
		v^2(r_v(\tau_2,\tau_1)-r_v(t_2,\tau_1))\Sigma_{I \backslash V,V} &
		v^2(r_v(\tau_2,\tau_2)-r_v(t_2,\tau_2))\Sigma_{I \backslash V, I \backslash V}
	\end{array}\right\Arrowvert \\
&\le& C_{M,\epsilon}+\epsilon(\tau_2-\tau_1).
\EQNY
Thus, we conclude that  for any small $\epsilon>0$
	there exist $C_{M,\epsilon}>0$ and $v_{M,\epsilon}>0$ such that for $v \geq v_{M,\epsilon}$,
inequality \kk{\eqref{eq:Pavel-1} \kls{with $H=1/2$} holds.}

\kk{\underline{ Inequality \eqref{eq:Pavel-2}}}.
For any $t_1 \in [\tau_1,\tau_1 +1]$ and $t_2 \in [\tau_2,\tau_2+1]$, we have, for the covariance matrix of $\vk{\chi}_{v}(t_1,t_2)$, that
\BQNY
\text{Cov} (\vk{\chi}_{v}(t_1,t_2)) )=v^2 \LT( \widetilde{R}_v(t_1,t_2,t_1,t_2) - \widetilde{R}_v(t_1,t_2,\tau_1,\tau_2)\widetilde{\Sigma}_v^{-1}(\tau_1,\tau_2) \widetilde{R}_v(\tau_1,\tau_2,t_1,t_2) \RT).
\EQNY
Using \eqref{r_uasymphbig} and some  calculations, we can show that,
for any small $\epsilon>0$ 
	there exist $C_{M,\epsilon}>0$ and $v_{M,\epsilon}>0$ such that for all $v \geq v_{M,\epsilon}$ and all $t_1 \in [\tau_1,\tau_1 +1]$, $t_2 \in [\tau_2,\tau_2+1]$,
\BQNY
|| \text{Cov} (\vk{\chi}_{v}(t_1,t_2)) ) || \le C_{M,\epsilon}+\epsilon(\tau_2-\tau_1).
\EQNY
\COM{
	it holds that
	\BQN
	\label{Pavel_var_2'}
	\mathrm{Var}(  \overline{ \vk w}_F^{\top} \vk{\chi}_{v,F}(t_1,t_2))
	&=& \overline{ \vk w}_{F \cap V}^{\top} \Sigma_{F \cap V,F \cap V}  \overline{ \vk w}_{F \cap V} r_v(t_1,t_1)
	+2  \overline{ \vk w}_{F \cap(I\setminus V)}^{\top} \Sigma_{F \cap (I\setminus V),F \cap V}  \overline{ \vk w}_{F \cap V} r_v(t_1,t_2)
	\notag
	\\
	&&
	+
	 \overline{ \vk w}_{F \cap (I\setminus V)}^{\top} \Sigma_{F \cap (I\setminus V),F \cap (I\setminus V)}  \overline{ \vk w}_{F \cap (I\setminus V)} r_v(t_2,t_2)
	\notag
	\\
	&\leq& C_{M,\epsilon}+\epsilon|\tau_2-\tau_1|.
	\EQN
}
Thus, for any small $\epsilon>0$
	there exist $C_{M,\epsilon}>0$ and $v_{M,\epsilon}>0$ such that for $v \geq v_{M,\epsilon}$
\kk{inequality \eqref{eq:Pavel-2} \kls{with $H=1/2$}  holds.}

\kk{\underline{ Inequality \eqref{eq:Pavel-3}}}.
We have, using \eqref{cond}, that, for any $F\subseteq I,$
\BQN
\Var\left(\overline{\vk{w}}_F^{\top} \chi_{v,F}(t_1,t_2)-\overline{\vk{w}}_F^{\top} \chi_{v,F}(s_1,s_2)\right)
\leq v^2 \Var \left(\overline{\vk{w}}_F^{\top} \widetilde{\vk{X}}_{v,F}(t_1,t_2)-\overline{\vk{w}}_F^{\top} \widetilde{\vk{X}}_{v,F}(s_1,s_2)\right).
\EQN
Further, 
it follows that there exist $C_M>0$, $v_M>0$ such that for all $v \geq v_M$ and all $t_1,s_1 \in [\tau_1,\tau_1+1]$, $t_2,s_2 \in [\tau_2,\tau_2 + 1]$, 
\BQN
\label{Pavel_Hoelder'}
&&
v^2 \Var \left(\overline{\vk{w}}_F^{\top} \widetilde{\vk{X}}_{v,F}(t_1,t_2)-\overline{\vk{w}}_F^{\top} \widetilde{\vk{X}}_{v,F}(s_1,s_2)\right)
\notag
\\
&&
\leq
2v^2 \mathrm{Var}
\left(\overline{\vk{w}}_{F\cap V}^{\top} {\vk{X}}_{v, F \cap V}(t_1)
-
\overline{\vk{w}}_{F\cap V}^{\top}
{\vk{X}}_{v, F \cap V}(s_1)\right)
\notag
\\
&&
\quad
+
2v^2 \mathrm{Var}
\left(\overline{\vk{w}}_{F \cap (I \backslash V)}^{\top} {\vk{X}}_{v, F \cap (I \backslash V)}(t_2)
-
\overline{\vk{w}}_{F \cap (I \backslash V)}^{\top}
{\vk{X}}_{v, F \cap (I \backslash V)}(s_2)\right)
\notag
\\
&&
\leq C_M(|t_1-s_1|+|t_2-s_2|),
\EQN
where the last inequality follows by an application of \eqref{r_uhoelderhbig}.
\kk{The above confirms that inequality \eqref{eq:Pavel-3}  \kls{with $H=1/2$} is \qq{satisfied}.
\\
To sum up,
we have checked that the conditions of \cite[Lemma 8]{Iev24} are satisfied, and therefore,}
for any small $\epsilon>0$ there exist $C_{M,\epsilon},v_{M,\epsilon}>0$ such that for all $v \geq v_{M,\epsilon}$,
\BQN
\int_{\R^{\abs{I}}} e^{(\widetilde{\vk{w}}_v(\tau_1,\tau_2))^{\top} \vk{x}}
\pk{\exists t_1 \in [\tau_1,\tau_1+1],t_2 \in [\tau_2,\tau_2+1]:\vk{\chi}_{v}(t_1,t_2)>\vk{x}}
d\vk{x} \leq e^{C_{M,\epsilon}+\epsilon(\tau_2-\tau_1)}.\label{eq:Int_w_2}
\EQN
It remains to estimate $\varphi_{\widetilde{\Sigma}_v(\tau_1,\tau_2)}\left(v \vk{b}_v(\tau_1,\tau_2)\right)$.  By \eqref{r_uasymphbig}, we have, for any small $\epsilon>0$ there exist $C_{M,\epsilon},v_{M,\epsilon}>0$ such that, for all $v \geq v_{M,\epsilon}$, 
\BQN\label{phi_bound}
\begin{aligned}
\varphi_{\widetilde{\Sigma}_v(\tau_1,\tau_2)}\left(
v \vk{b}_v(\tau_1,\tau_2)\right)
&
\leq C_{M,\epsilon}\exp\left(-\frac{1}{2}
v^2 \vk{b}_v^{\top}(\tau_1,\tau_2)
\widetilde{\Sigma}_v^{-1}(\tau_1,\tau_2)
\vk{b}_v(\tau_1,\tau_2)
\right)
\\&
\leq
C_{M,\epsilon} \exp\left(-\frac{1}{2}
v^2 \vk{b}_v^{\top}(\tau_1,\tau_2)
\Sigma^{-1}\vk{b}_v(\tau_1,\tau_2)
+\epsilon(\tau_2-\tau_1)
\right).
\end{aligned}
\EQN
On the other hand, we have by Lagrange's mean value theorem that, for some $x \in [\tau_1,\tau_2]$,  
\BQNY
\vk{b}_v^{\top}(\tau_1,\tau_2)\Sigma^{-1} \vk{b}_v(\tau_1,\tau_2)=
\vk{b}_v^{\top}(\tau_1,\tau_1)\Sigma^{-1} \vk{b}_v(\tau_1,\tau_1)
+2v^{-2}(\tau_2-\tau_1) \vk{b}_v^{\top}(\tau_1,x)\Sigma^{-1}
\begin{pmatrix}
	\vk{0}_V\\
	\vk{\zeta}'_{I \backslash V}(t_0+xv^{-2})
\end{pmatrix}.
\EQNY
Therefore, for any small $\epsilon>0$, there exists $v_{M,\epsilon}$ such that, for all $v \geq v_{M,\epsilon}$
\BQN
\label{quad_form_asymp}
v^2
\abs{
\vk{b}_v^{\top}(\tau_1,\tau_2)\Sigma^{-1} \vk{b}_v(\tau_1,\tau_2)-
g_I(t_0+\tau_1 v^{-2})
-2v^{-2}(\tau_2-\tau_1)\vk{\zeta}_I(t_0)^{\top}\Sigma_{II}^{-1}
\begin{pmatrix}
	\vk{0}_V\\
	\vk{\zeta}'_{I \backslash V}(t_0)
\end{pmatrix}
}
\leq \epsilon(\tau_2-\tau_1).
\EQN
Now take $V = \{i \in I: \zeta'_i(t_0) \leq 0\}$. Since $\vk{\zeta}'_I(t_0) \neq \vk{0}_I$ \ksn{(due to $H\vk{\nu}_I \neq (1-H)t_0\vk{\mu}_I$)}, $\Sigma_{II}^{-1} \vk{\zeta}_I(t_0)> \vk{0}_I$  and $2\vk{\zeta}_I(t_0)^{\top}\Sigma_{II}^{-1}\vk{\zeta}'_I(t_0)=g'_I(t_0)=0$, we know that both $V$ and $I \backslash V$ are non-empty. Hence, using once again  that $\Sigma_{II}^{-1} \vk{\zeta}_I(t_0)> \vk{0}_I$, we obtain
\BQN
\label{positive_form}
\vk{\zeta}_I(t_0)^{\top} \Sigma_{II}^{-1}
\begin{pmatrix}
	\vk{0}_V\\
	\vk{\zeta}'_{I \backslash V}(t_0)
\end{pmatrix}
>0.
\EQN
Therefore, combining \eqref{phi_bound}-\eqref{positive_form} we have that, for any small  $\epsilon>0$ there exist $C, C_{M,\epsilon},v_{M,\epsilon}>0$ such that for all $v \geq v_{M,\epsilon}$,
\BQNY
\varphi_{\widetilde{\Sigma}_v(\tau_1,\tau_2)}(v\vk{b}_v(\tau_1,\tau_2))
&\leq& C_{M,\epsilon} \exp\left(-\frac{v^2}{2}g_I(t_0+\tau_1 v^{-2})- C(\tau_2-\tau_1)\right)
\\
&\leq&
C_{M,\epsilon}
\exp\left(-\frac{v^2}{2}g_I(t_0)- C(\tau_2-\tau_1)\right).
\EQNY
This, combined with \ccP{\eqref{int} and} \eqref{eq:Int_w_2}, establishes the claim of the lemma. 
\QED

\section*{Declarations}

{\bf Acknowledgement:}	We are thankful to the reviewer for his/her comments, which have helped to improve the presentation of the manuscript.

{\bf Authors' contributions.} L.J. and S.N. wrote the main manuscript text based on several discussions among all authors and K.D. prepared the Introduction. All authors reviewed and revised the manuscript.

{\bf Funding.} Financial support from the Swiss National Science Foundation Grant 200021-196888 is kindly acknowledged. K. D\c{e}bicki was partially supported by NCN Grant No 2018/31/B/ST1/00370 (2019-2024).

{\bf Data availability.}	 Data sharing not applicable to this article as no data sets were generated or analysed during	the current study.

{\bf Ethical approval.} Not applicable. 

{\bf Conflict of interest.}  The authors declare that they have no conflicts of interest to this work.


\bibliographystyle{plain} 

\bibliography{EEEA.bib}       

\end{document}